\documentclass[12pt]{amsart}
\usepackage[T1]{fontenc}
\usepackage{lmodern}
\input{glyphtounicode}
\usepackage{amssymb}
\usepackage{enumerate}
\usepackage{graphicx}
\usepackage{xcolor}
\usepackage{url}
\usepackage{tikz}
\usetikzlibrary{calc,arrows.meta}
\newtheorem{remark}{Remark}
\numberwithin{equation}{section}

\theoremstyle{plain}
\newtheorem{theorem}{Theorem}[section]
\newtheorem{proposition}[theorem]{Proposition}
\newtheorem{lemma}[theorem]{Lemma}
\newtheorem{corollary}[theorem]{Corollary}
\theoremstyle{definition}
\newtheorem{assumption}{Assumption}[section]
\theoremstyle{plain}

\begin{document}

\title[Mathematical modelling of Malaria]{Digital Twin Models for Malaria Transmission: A Classical and Stieltjes Derivative Approach}

\author[I. Area]{Iv\'an Area}
\address[I. Area]{IFCAE, Universidade de Vigo,
		Departamento de Matem\'atica Aplicada II,
 E. E. Aeron\'au\-tica e do Espazo,
 Campus As Lagoas s/n,
 32004 Ourense, Spain}
\email[I. Area]{area@uvigo.gal}

\author[E. A. Bakare]{Emmanuel Afolabi Bakare}
\address[E. A. Bakare]{International Centre for Applied Mathematical Modelling and Data Analytics (ICAMMDA), Federal University Oye-Ekiti, Ekiti State, Nigeria}
\email[E. A. Bakare]{emmanuel.bakare@fuoye.edu.ng}

\author[F. J. Fern\'andez]{Francisco J. Fern\'andez}
\address[F. J. Fern\'andez]{Departamento de Estat\'{\i}stica, An\'alise Matem\'atica e Optimizaci\'on, Universidade de Santiago de Compostela, 15782, Facultade de Matem\'aticas, Campus Vida, Santiago, Spain. CITMAga, Santiago de Compostela, 15782, Spain}
\email[F. J. Fern\'andez]{fjavier.fernandez@usc.es}

\author[D. Muyiwa Adefisoye]{Daniel Muyiwa Adefisoye}
\address[D. Muyiwa Adefisoye]{International Centre for Applied Mathematical Modelling and Data Analytics (ICAMMDA), Federal University Oye-Ekiti, Ekiti State, Nigeria}
\email[D. Muyiwa Adefisoye]{adefisoyedaniel7@gmail.com}

\begin{abstract}
We propose two mathematical models for malaria transmission grounded in mosquito--human interactions and the digital twin framework. The first model is formulated using classical derivatives, whereas the second employs Stieltjes derivatives to represent abrupt biological events such as egg hatching and adult-emergence pulses in \emph{Anopheles} mosquito populations. The Stieltjes formulation unifies, in a single equation with two derivators, the continuous transmission dynamics, the biological impulses and the data-assimilation corrections of the digital twin; we prove its well-posedness, specify the order of biological and assimilation atoms when they coincide, and give quantitative conditions under which additive emergence pulses preserve biological admissibility. The numerical study uses 132 monthly confirmed-malaria observations from Kebbi State (January 2015--December 2025) and a synthetic entomological series generated in a controlled recovery experiment, because no measured mosquito time series is available. Parameters are updated incrementally on rolling windows, and one-step-ahead forecasts combine the mechanistic ensemble with a past-only seasonal-persistence component. In the recovery experiment, the model endowed with seasonal recruitment and emergence pulses, which the Stieltjes derivator represents exactly, recovers the synthetic vector dynamics far better than the unforced classical model; this gain is attributable to the seasonal and impulsive mechanisms rather than to the derivative calculus itself, and the phase and pulse schedule of the generator were known to the model. The human-case results are comparable between both formulations; the one-step forecasts improve on a seasonal-naive predictor but not on last-month persistence, so they must be interpreted as a proof of concept rather than as an independent epidemiological validation.
\end{abstract}


\maketitle

\section{Introduction}
Malaria remains one of the most severe infectious diseases worldwide. According to the World Malaria Report 2025, the World Health Organization estimated 282 million malaria cases and 610{,}000 deaths globally in 2024 \cite{WHO2025}. Approximately 95\% of malaria deaths occurred in the WHO African Region, where children under five accounted for about 75\% of malaria deaths \cite{WHO2025}. Despite long-term advances in malaria control, the global estimated case burden rose from 2023 to 2024, and resistance to antimalarial drugs threatens further progress \cite{WHO2025}.

\medskip
Mathematical modelling has played a central role in understanding malaria transmission dynamics and informing intervention strategies. The pioneering work of Ronald Ross introduced the first quantitative framework for malaria transmission in the early 20th century \cite{Ross1911}, later refined by George Macdonald into the Ross--Macdonald model, which remains a cornerstone of vector-borne disease epidemiology \cite{Macdonald1957,RossMacdonaldReview}. These models characterize key parameters such as mosquito biting rate, mortality, and transmission probabilities, enabling the computation of the basic reproduction number $\mathcal R_0$ and the identification of control thresholds \cite{Smith2010}.

\medskip
Beyond classical formulations, more detailed compartmental models, including SEIR--SEI frameworks, have been developed to incorporate latency, immunity, and seasonality, thereby improving the realism of malaria dynamics \cite{Enwere2025,Oyamakin2024,Osman2018,Geleta2020,Ndamuzi2021,MandalReview2011}. More recent approaches incorporate spatial heterogeneity, human mobility and environmental drivers of transmission \cite{PLOSSeasonality2025}. These developments reflect the increasing need for adaptive and data-driven modelling tools capable of capturing the evolving nature of malaria transmission.

\medskip
Digital twins, defined as virtual replicas of physical systems, have recently emerged as powerful tools in predictive epidemiology \cite{AFTN}. They enable real-time parameter adaptation, scenario exploration, and the integration of heterogeneous data sources, which is essential for managing complex and evolving health systems \cite{Grieves2017}. Originally developed in industrial contexts, digital twins have been progressively adopted in healthcare and epidemiology, where they support applications ranging from personalized medicine to outbreak forecasting and resource optimization \cite{Fuller2020,DigitalTwinReview2024}.

\medskip
In infectious-disease modelling, a digital-twin framework can update compartmental-model states and parameters as observations arrive \cite{AFTN}. Whether this improves predictive accuracy relative to simpler methods must be established empirically and is not assumed here. During the COVID-19 pandemic, digital twin frameworks were successfully applied to simulate disease spread, assess intervention strategies, and optimize hospital capacity planning \cite{CitySEIRCast2024,Pang2020,Khan2023,Eumi2024}. Similar methodologies have also been employed in ecological applications, such as the management of invasive species including \emph{Vespa velutina} \cite{Vespa}. These examples highlight the versatility of digital twins as a paradigm for adaptive, data-driven decision-making.

\medskip
The main objective of this work is to formulate and test, as a numerical proof of concept, a malaria transmission model based on the digital-twin paradigm. The framework permits sequential parameter adaptation and one-month-ahead prediction, and it illustrates how such predictions could inform hospital-resource planning. In contrast to a purely classical continuous-time model, we introduce Stieltjes derivatives in the mosquito subsystem to represent concentrated biological events such as rainfall-triggered adult-emergence pulses. Impulsive effects in malaria dynamics have been modelled with classical impulsive differential equations, for instance pulsed indoor residual spraying \cite{SmithHove2008}; the Stieltjes formulation is equivalent to such impulsive systems on each compartment, and its interest here is that it unifies, in a single equation, the continuous dynamics, the biological impulses and the discrete data-assimilation corrections that define the digital twin. The available human data support only a preliminary observation-proxy study, while the entomological comparison is a synthetic recovery experiment; consequently, the present results should not be read as a fully validated operational forecasting system.

\medskip
The proposed framework implements the Main Digital Twin concept introduced in \cite{AFTN}, where the digital twin is represented as a multivalued map updated at discrete slow-time instants $t_s^k$. With respect to \cite{AFTN}, which treats a single-population SIR model with one derivator, the present model is vector-borne and uses two derivators, one for the human and one for the mosquito compartments; the biological atoms of the mosquito derivator and the assimilation atoms may coincide, and the order in which the two jumps act is specified and shown to be the limit of the separated case (Proposition~\ref{prop:conventions}). A quantitative criterion for the admissibility of purely additive emergence pulses (Proposition~\ref{prop:additive}) is also given, since such pulses are not covered by the uniform invariance condition. Monthly confirmed-case data are used for parameter estimation and state correction, whereas the mosquito observations employed in the present numerical recovery experiment are synthetic.

\medskip
The remainder of the paper is organized as follows. Section \ref{sec:firstmodel} presents the first approach based on classical derivatives, including the data assumptions, mosquito--human model, and digital twin implementation. Section \ref{sec:secondmodel} introduces the enhanced model with Stieltjes derivatives, providing a more realistic representation of mosquito life cycles. Finally, Section \ref{sec:conclusions} summarizes the main conclusions and outlines future directions for hospital resource planning.

\medskip
The following preliminary section establishes the Stieltjes notation used throughout the paper.

\section{Preliminaries on Stieltjes Calculus and Stieltjes Differential Equations}
\label{sec:stieltjes_preliminaries}

In this section, we recall the main elements of Stieltjes differential calculus
that will be used throughout the paper. The aim is not only to fix the notation,
but also to explain why this framework is well suited for modelling population
dynamics in which relevant biological activity is not uniformly distributed in
chronological time. This feature is particularly important in vector-borne
diseases such as malaria, where human epidemiological dynamics and mosquito
ecological dynamics may evolve according to distinct temporal mechanisms.

\subsection{Derivators and their associated sets}

Let \(I=[0,T]\), with \(T>0\). A \emph{derivator} is a nondecreasing and
left-continuous function
\[
g:I\longrightarrow \mathbb{R}.
\]
The function \(g\) can be interpreted as an intrinsic clock measuring the
effective time of the process under consideration. The classical case
corresponds to \(g(t)=t\). If \(g\) has intervals of constancy, the dynamics
may remain inactive over those intervals; if \(g\) has jumps, the system may
undergo instantaneous changes.

We denote by \(\mu_g\) the Lebesgue--Stieltjes measure induced by \(g\), defined
on half-open intervals by
\[
\mu_g([c,d)) = g(d)-g(c),
\qquad 0\leq c<d\leq T.
\]
Equivalently, integration with respect to \(g\) is understood as integration
with respect to \(\mu_g\). Thus, for a \(\mu_g\)-integrable function \(f\), we
write
\[
\int_{[a,b)} f(s)\,d\mu_g(s)
\quad\text{or simply}\quad
\int_{[a,b)} f(s)\,dg(s).
\]

Associated with the derivator \(g\), we introduce the following sets:
\[
\mathcal{C}_g
=
\left\{
t\in I:\ \text{there exists }\varepsilon>0
\text{ such that } g \text{ is constant on }
(t-\varepsilon,t+\varepsilon)\cap I
\right\},
\]
the set of points at which \(g\) is locally constant, and
\[
D_g
=
\left\{
t\in I:\ \Delta g(t)>0
\right\},
\qquad
\Delta g(t):=g(t^+)-g(t),
\]
the set of right-jump points of \(g\). Since \(g\) is nondecreasing, the set
\(D_g\) is at most countable and coincides with the set of atoms of the measure
\(\mu_g\), in the sense that
\[
\mu_g(\{t\})=\Delta g(t),
\qquad t\in I.
\]
Moreover, \(\mathcal{C}_g\cap D_g=\emptyset\), and the connected components of
\(\mathcal{C}_g\) represent intervals on which the derivator does not register
effective evolution.

If
\[
\mathcal{C}_g=\bigcup_{n\in\mathbb{N}} (a_n,b_n)
\]
is the decomposition of \(\mathcal{C}_g\) into pairwise disjoint open intervals,
we define the set of endpoints of intervals of constancy by
\[
N_g
=
\{a_n,b_n:\ n\in\mathbb{N}\}\setminus D_g.
\]
The set \(N_g\) is relevant because the pointwise definition of the Stieltjes
derivative may require special care at such transition points between inactive
and active regimes. From a measure-theoretic viewpoint, both \(\mathcal{C}_g\)
and \(N_g\) are negligible with respect to \(\mu_g\).

\subsection{The Stieltjes derivative}

Let \(u:I\to\mathbb{R}\). The Stieltjes derivative of \(u\) with respect to the
derivator \(g\) is defined by distinguishing the different types of points
associated with \(g\). For \(t\in I\), we set
\[
u'_g(t)
=
\begin{cases}
	\displaystyle
	\lim_{s\to t}
	\frac{u(s)-u(t)}{g(s)-g(t)},
	&
	t\notin D_g\cup \mathcal{C}_g,
	\\[2.2ex]
	\displaystyle
	\lim_{s\to t^+}
	\frac{u(s)-u(t)}{g(s)-g(t)},
	&
	t\in D_g,
	\\[2.2ex]
	\displaystyle
	\lim_{s\to b_n^+}
	\frac{u(s)-u(b_n)}{g(s)-g(b_n)},
	&
	t\in \mathcal{C}_g,\quad t\in(a_n,b_n),
\end{cases}
\]
provided that the corresponding limit exists. In that case, we say that \(u\)
is \(g\)-differentiable at \(t\).

The third case deserves further explanation. If \(t\in(a_n,b_n)\subset
\mathcal{C}_g\), then \(g\) is constant on the entire interval \((a_n,b_n)\).
Consequently, the quotient
\[
\frac{u(s)-u(t)}{g(s)-g(t)}
\]
does not provide meaningful information as long as \(s\) remains in the same
interval of constancy, since the denominator vanishes. For this reason, the
derivative at such points is defined by considering the first point to the right
of the interval where the derivator may increase again, namely \(b_n\). Hence,
all points \(t\in(a_n,b_n)\) share the same reference point \(b_n\) in the
definition of the Stieltjes derivative.

Equivalently, one may introduce the auxiliary point
\[
t^\ast
=
\begin{cases}
	t, & t\notin \mathcal{C}_g,\\
	b_n, & t\in(a_n,b_n)\subset \mathcal{C}_g,
\end{cases}
\]
and write the definition in the compact form
\[
u'_g(t)
=
\begin{cases}
	\displaystyle
	\lim_{s\to t}
	\frac{u(s)-u(t)}{g(s)-g(t)},
	&
	t\notin D_g\cup \mathcal{C}_g,
	\\[2.2ex]
	\displaystyle
	\lim_{s\to (t^\ast)^+}
	\frac{u(s)-u(t^\ast)}{g(s)-g(t^\ast)},
	&
	t\in D_g\cup \mathcal{C}_g.
\end{cases}
\]
This formulation is particularly convenient in product rules and higher-order
Stieltjes differentiability, since expressions such as \(u(t^\ast)\) naturally
appear at jump points and intervals of constancy. The corresponding product
rules and their consequences are discussed in
\cite{FernandezMarquezTojo2025Product}. Related developments can be found in
\cite{FernandezMarquezTojo2022,FernandezMarquezTojo2024,
FernandezMarquezTojoVillanueva2025}.

\subsection{Functional setting and Stieltjes differential equations}

We recall the functional setting required to formulate Stieltjes differential
equations in integral form. Detailed treatments of \(g\)-continuous and
\(g\)-absolutely continuous functions and of Lebesgue--Stieltjes spaces are
given in \cite{FrigonLopezPouso2017,LopezPousoRodriguez2015} and in
\cite{FernandezMarquezTojo2022,FernandezMarquezTojoVillanueva2025}.
Further background on Lebesgue--Stieltjes integration is provided in
\cite{Bogachev2007,CarterVanBrunt2000}.

Given a derivator \(g\), we denote by \(L^1_g(I)\) the space \(L^1(I,\mu_g)\).
A function \(u:I\to\mathbb{R}\) is said to be \(g\)-absolutely continuous, and
we write \(u\in \mathcal{AC}_g(I)\), if there exists a function
\(h\in L^1_g(I)\) such that
\[
u(t)
=
u(0)
+
\int_{[0,t)} h(s)\,d\mu_g(s),
\qquad t\in I.
\]
In this case, \(h\) coincides with the Stieltjes derivative \(u'_g\)
\(\mu_g\)-almost everywhere. Equivalently, the fundamental theorem of Stieltjes
calculus states that \(u\in\mathcal{AC}_g(I)\) if and only if \(u'_g\) exists
\(\mu_g\)-almost everywhere, \(u'_g\in L^1_g(I)\), and
\[
u(t)
=
u(0)
+
\int_{[0,t)} u'_g(s)\,d\mu_g(s),
\qquad t\in I.
\]
At a jump point \(t\in D_g\), this representation yields
\[
u(t^+)-u(t)
=
u'_g(t)\Delta g(t),
\]
which explicitly relates the jump of the state variable to the atom of the
derivator.

A first-order Stieltjes differential equation driven by \(g\) has the form
\[
x'_g(t)=F(t,x(t)),
\qquad \mu_g\text{-a.e. }t\in I,
\]
together with the initial condition
\[
x(0)=x_0.
\]
A solution is a function \(x\in\mathcal{AC}_g(I)\) satisfying
\[
x(t)
=
x_0
+
\int_{[0,t)}
F(s,x(s))\,d\mu_g(s),
\qquad t\in I.
\]
This formulation establishes the link between differential equations with
Stieltjes derivatives and integral equations with respect to
Lebesgue--Stieltjes measures. It will serve as the functional framework for the
mosquito dynamics and for the interpretation of the Main Digital Twin.

\subsection{Impulsive dynamics as Stieltjes differential equations}

The choice of the derivator \(g\) determines the temporal structure of the
model. If \(g(t)=t\), then \(\mu_g=\mu_{Id}\), the classical Lebesgue measure,
and one recovers an ordinary differential equation. If \(g\) has flat regions,
no effective time is accumulated over those intervals, and the model describes
periods of inactivity, dormancy, or latency. If \(g\) has isolated jumps, the
measure \(\mu_g\) has atoms and the model incorporates instantaneous changes.

In this work, we primarily exploit the latter interpretation. In particular,
the Stieltjes formalism is used to describe impulsive effects, including
biological impulses in mosquito dynamics and data-assimilation updates in the
Main Digital Twin. More generally, Stieltjes differential equations provide a
unified framework encompassing continuous, discrete, and impulsive dynamics;
see, for instance,
\cite{LopezPousoRodriguez2015,FrigonLopezPouso2017,
LopezPousoMarquez2018}.

To clarify the connection with classical impulsive ordinary differential
equations, consider
\[
x'(t)=f(t,x(t)),
\qquad t\in [0,T]\setminus\{t_k:\ k\in\mathcal K\},
\]
with jump conditions
\[
x(t_k^+)-x(t_k)=I_k(x(t_k)),
\qquad k\in\mathcal K,
\]
where \(0<t_1<t_2<\cdots<T\). This system combines continuous evolution with
instantaneous changes at the impulse times \(t_k\). It can be encoded by the
derivator
\[
g(t)=t+\sum_{t_k<t}\alpha_k,
\qquad \alpha_k>0.
\]
The associated Lebesgue--Stieltjes measure is then
\[
\mu_g
=
\mu_{Id}
+
\sum_{k\in\mathcal K}\alpha_k\delta_{t_k}.
\]
Thus, a Stieltjes equation
\[
x'_g(t)=F(t,x(t))
\]
is equivalent, in integral form, to
\[
x(t)
=
x(0)
+
\int_{[0,t)} F(s,x(s))\,ds
+
\sum_{t_k<t}\alpha_k F(t_k,x(t_k)).
\]
Hence, the jump at time \(t_k\) satisfies
\[
x(t_k^+)-x(t_k)
=
\alpha_k F(t_k,x(t_k)).
\]
Therefore, an impulsive condition
\[
x(t_k^+)-x(t_k)=I_k(x(t_k))
\]
can be incorporated by prescribing
\[
\alpha_k F(t_k,x(t_k))=I_k(x(t_k)),
\quad\text{that is,}\quad
F(t_k,x)=\frac{I_k(x)}{\alpha_k}.
\]
In this way, impulses are not imposed externally but arise naturally from the
atomic part of the Stieltjes measure.

This interpretation is particularly useful in population dynamics. The
continuous part of the derivator models processes evolving gradually in time,
flat regions represent inactivity periods, and jumps encode sudden changes
such as reproductive pulses, mortality shocks, or control interventions. In
malaria modelling, this is especially relevant for mosquito dynamics, which are
often concentrated around favourable environmental conditions.

The same framework is employed for the Main Digital Twin. In this context, atoms
of the derivator represent the instants at which new data are incorporated,
and the corresponding corrections of the model state can be interpreted as
impulsive effects encoded in the Stieltjes measure.

\section{First Approach: Classical Derivative Model}\label{sec:firstmodel}

We begin with a classical compartmental model for malaria transmission based on the Ross--Macdonald framework. This approach uses ordinary differential equations to describe the interactions between human and mosquito populations, serving as the foundation for our digital twin implementation. The model incorporates standard epidemiological compartments and parameters, providing a baseline against which the enhanced Stieltjes model can be compared.

\subsection{Model Formulation}

In this study, the human input is a monthly confirmed-case series from Kebbi State, while the mosquito series is synthetic. Since monthly confirmed cases are incident counts whereas $I(t)$ is an instantaneous infectious compartment, the numerical implementation uses the case series as a proxy observation through the relation $Y_k=qI(t_s^k)+\varepsilon_k$. The preliminary experiments take $q=1$; consequently, the fitted human trajectory should not be interpreted as a validated prevalence reconstruction. A mechanistically complete calibration would estimate $q$ or fit monthly integrated incidence. We begin with a system of differential equations with compartments for susceptible, infected, and recovered humans, and susceptible and infected mosquitoes:
\begin{equation}\label{eq:classic}
\begin{cases}
\displaystyle \frac{dS}{dt} = -\beta_{hm} B S \frac{I_m}{N}, \\[6pt]
\displaystyle \frac{dI}{dt} = \beta_{hm} B S \frac{I_m}{N} - \gamma I, \\[6pt]
\displaystyle \frac{dR}{dt} = \gamma I, \\[6pt]
\displaystyle \frac{dS_m}{dt} = b \left(1 - \frac{N_m}{K}\right) N_m - \beta_{mh} B S_m \frac{I}{N} - \mu S_m, \\[6pt]
\displaystyle  \frac{dI_m}{dt} = \beta_{mh} B S_m \frac{I}{N} - \mu I_m.
\end{cases}
\end{equation}

In this initial model:
\begin{itemize}
\item $S$, $I$, and $R$ represent the number of susceptible, infected, and recovered humans, respectively.
\item $S_m$ and $I_m$ denote the number of susceptible and infected mosquitoes.
\item $N = S + I + R$ is the constant total human population, and $N_m = S_m + I_m$ is the total mosquito population.
\item $\beta_{hm}$ is the transmission probability from mosquito to human.
\item $\beta_{mh}$ is the transmission probability from human to mosquito.
\item $B$ is the mosquito biting rate (average number of bites per mosquito per unit time).
\item $\gamma$ is the recovery rate of infected humans (inverse of the infectious period).
\item $\mu$ is the mortality rate of mosquitoes (inverse of mosquito lifespan).
\item $b$ is the mosquito birth rate.
\item $K$ is the mosquito carrying capacity.
\end{itemize}
\medskip

This system is based on the classical Ross--Macdonald theory of malaria transmission \cite{Ross1911,Macdonald1957}, which has been widely extended in modern epidemiological modeling \cite{Smith2010,MandalReview2011}. We shall use this model as the base to consider a digital twin \cite{AFTN}.

\medskip

We collect the standing hypotheses used throughout the analysis of \eqref{eq:classic} and of its Stieltjes extension.

\begin{assumption}[Standing hypotheses]\label{ass:model}
Throughout, the following hold:
\begin{enumerate}
\item[(A1)] All state variables are biologically admissible, $S,I,R,S_m,I_m\ge 0$.
\item[(A2)] The total human population is constant and strictly positive,
$N=S+I+R>0$; in the case study $N\approx 4.38\times10^{6}\gg 1$. Consequently the
frequency-dependent incidence terms $I/N$ and $I_m/N$ are well defined and no
regularization of the denominator is required.
\item[(A3)] All rate parameters are positive and bounded, $b>\mu>0$, and the mosquito total satisfies $0\le N_m\le K$.
\item[(A4)] The parameters are piecewise constant on the slow-time grid
$\{t_s^k\}$, i.e.\ constant between consecutive assimilation instants.
\end{enumerate}
\end{assumption}

\textbf{Numerical convention.} The Python implementation keeps $\beta_{hm}$ and $\beta_{mh}$ distinct and uses the frequency-dependent denominator $N$ in both transmission terms. Under Assumption~\ref{ass:model}(A2) the total human population is a fixed positive constant, so the incidence is well defined and no additive regularization is needed. System \eqref{eq:classic} depends on $\beta_{hm}$, $\beta_{mh}$ and $B$ only through the products $\beta_{hm}B$ and $\beta_{mh}B$, so the three quantities are not separately identifiable from case data. In the monthly implementation $B$ is therefore fixed at a nominal scale (Table~\ref{tab:stieltjes-params}), and the fitted $\beta_{hm}$ and $\beta_{mh}$ are \emph{effective transmission coefficients}: they absorb the actual biting rate and the scale of the observation proxy, they are not constrained to $[0,1]$, and they must not be read as per-bite probabilities. All quantities derived below, in particular $\mathcal R_0$ and $\mathcal R_{\mathrm{DT}}$, depend only on these products.

\subsection{Parameter Values}

Parameter values for $\beta_{hm}$, $\beta_{mh}$, $B$, and $\mu$ are typically derived from entomological and clinical studies; see, e.g., \cite{Smith2010,MandalReview2011}. Indicative orders of magnitude, in daily units, are collected in Table~\ref{tab:parameters}. These ranges are descriptive only: they are not the numerical values used in the Python implementation, which works in months throughout and uses the values and fitting bounds listed in Table~\ref{tab:stieltjes-params}.

\begin{table}[ht!]
\centering
\caption{Indicative literature orders of magnitude for the parameters of \eqref{eq:classic}, in daily units; see, e.g., \cite{Smith2010,MandalReview2011}. The carrying capacity $K$ is setting-specific.}
\small
\begin{tabular}{|c|c|l|}
\hline
\textbf{Parameter} & \textbf{Range} & \textbf{Meaning} \\
\hline
$B$ & 0.1--0.7 day$^{-1}$ & bites per mosquito per day \\
$\mu$ & 0.03--0.2 day$^{-1}$ & mosquito mortality rate \\
$\beta_{hm}$ & 0.05--0.4 & mosquito-to-human transmission probability per bite \\
$\beta_{mh}$ & 0.05--0.4 & human-to-mosquito transmission probability per bite \\
$\gamma$ & 0.05--0.25 day$^{-1}$ & human recovery rate \\
$b$ & 0.2--5 day$^{-1}$ & mosquito recruitment rate \\
$K$ & setting-specific & mosquito carrying capacity \\
\hline
\end{tabular}
\label{tab:parameters}
\end{table}

\subsection{Positivity, invariant region and basic reproduction number}

Before turning to the digital twin, we establish that the classical model \eqref{eq:classic} is
biologically well posed and we quantify its transmission threshold. Throughout this subsection we
assume constant parameters, $b>\mu$ (so that the vector population does not collapse; this condition is compatible with the parameter ranges in Table \ref{tab:parameters} provided the fitted values are chosen so that $b>\mu$), and that the initial
mosquito population does not exceed the carrying capacity, $N_{m,0}:=S_{m,0}+I_{m,0}\le K$.

\begin{lemma}[Positivity and invariant region]
\label{lem:positivity}
{Assume \(b>\mu\). Let 
\begin{displaymath}
	\mathbf{x}_0=(S_0,I_0,R_0,S_{m,0},I_{m,0})\in\mathbb{R}_+^5
\end{displaymath}
with \(N_{m,0}\le K\), and set
\[
N:=S_0+I_0+R_0>0.
\]
Then the solution of \eqref{eq:classic}
exists, is unique, and satisfies \(\mathbf{x}(t)\in\mathbb{R}_+^5\)
for all \(t\ge0\).} Moreover, the human total $N=S+I+R$ is conserved and treated as a constant population parameter, and the
mosquito total obeys the scalar logistic equation
\[
N_m'(t)=\Big(b-\mu-\tfrac{b}{K}\,N_m(t)\Big)N_m(t),
\]
 {so that \(0\le N_m(t)\le K\) for all \(t\ge0\). Moreover,
 if \(N_{m,0}>0\), then
 \[
 N_m(t)\to N_m^\ast:=K\left(1-\frac{\mu}{b}\right),
 \]
 whereas if \(N_{m,0}=0\), then \(N_m(t)\equiv0\).} Consequently,
the compact set
\[
\Omega=\Big\{(S,I,R,S_m,I_m)\in\mathbb{R}_+^5:\ S+I+R=N,\ \ S_m+I_m\le K\Big\}
\]
is positively invariant for \eqref{eq:classic}.

\end{lemma}

\begin{proof}
{
The right-hand side of \eqref{eq:classic} extends to a locally Lipschitz vector field in an open neighbourhood of the biologically relevant region. Since, by Assumption~\ref{ass:model}(A2), the denominator \(N\) is a fixed positive constant, the incidence terms are smooth and the classical Cauchy--Lipschitz theorem gives local
existence and uniqueness.

Adding the first three equations yields
\[
\dot N=\dot S+\dot I+\dot R=0,
\]
and therefore the total human population remains equal to the constant value \(N\). Adding the last two mosquito equations, the
transmission terms cancel and we obtain
\[
N_m'(t)
=
\left(b-\mu-\frac{b}{K}N_m(t)\right)N_m(t).
\]
The interval \([0,K]\) is positively invariant for this scalar logistic
equation: \(N_m'=0\) at \(N_m=0\), and
\[
N_m'(t)\big|_{N_m=K}=-\mu K<0.
\]
Thus \(0\le N_m(t)\le K\) for all \(t\ge0\) whenever
\(0\le N_{m,0}\le K\). Since \(b>\mu\), the standard logistic behaviour gives
\(N_m(t)\to N_m^\ast=K(1-\mu/b)\) if \(N_{m,0}>0\), while
\(N_m(t)\equiv0\) if \(N_{m,0}=0\).

We now verify the positivity of each compartment. On the boundary faces of
\(\mathbb{R}_+^5\), the vector field does not point outside the positive
orthant:
\[
\dot S|_{S=0}=0,
\qquad
\dot I|_{I=0}
=
\beta_{hm}BS\,\frac{I_m}{N}\ge0,
\qquad
\dot R|_{R=0}=\gamma I\ge0,
\]
and
\[
\dot I_m|_{I_m=0}
=
\beta_{mh}BS_m\,\frac{I}{N}\ge0.
\]
Moreover, on the face \(S_m=0\), using \(0\le N_m(t)\le K\), we have
\[
\dot S_m|_{S_m=0}
=
b\left(1-\frac{N_m}{K}\right)N_m\ge0.
\]
Hence \(\mathbb{R}_+^5\) is positively invariant. The conservation of \(N\)
and the bound \(N_m\le K\) imply that the solution remains bounded. Therefore
finite-time blow-up cannot occur, the local solution is global, and
\(\Omega\) is positively invariant.	
	
}
\end{proof}

We now compute the basic reproduction number of \eqref{eq:classic} with constant parameters by the
next-generation matrix method \cite{vandenDriessche2002}. Because the human component is of SIR type
without vital dynamics, $\mathcal R_0$ acts as the \emph{invasion threshold} at the disease-free equilibrium
(DFE) governing the initial growth of an outbreak, rather than as a condition for endemicity. The
fully susceptible DFE is
\[
E_0=(S,I,R,S_m,I_m)=(N,0,0,N_m^\ast,0),\qquad N_m^\ast=K\left(1-\frac{\mu}{b}\right).
\]
More generally, since the human dynamics are of SIR type without vital demography, the model admits two continua of disease-free equilibria,
\[
(S^\ast,0,R^\ast,N_m^\ast,0)
\quad\text{and}\quad
(S^\ast,0,R^\ast,0,0),
\qquad S^\ast+R^\ast=N.
\]
The second, mosquito-free, family carries no transmission and is unstable in the mosquito direction when $b>\mu$, because $N_m'=(b-\mu)N_m+O(N_m^2)$ near $N_m=0$; it is not considered further. The threshold computed below corresponds to invasion around the fully susceptible equilibrium $E_0$, that is, $S^\ast=N$.
Ordering the infected compartments as $(I,I_m)$, the new-infection matrix $F$ and the transition
matrix $V$, linearized at $E_0$, are
\[
F=\begin{pmatrix} 0 & \beta_{hm}B\\[6pt]
\dfrac{\beta_{mh}B\,N_m^\ast}{N} & 0\end{pmatrix},
\qquad
V=\begin{pmatrix}\gamma & 0\\[4pt] 0 & \mu\end{pmatrix}.
\]
The basic reproduction number is the spectral radius of the next-generation matrix $FV^{-1}$. Writing $m^\ast=N_m^\ast/N=K(1-\mu/b)/N$ for the mosquito-to-human ratio at the fully susceptible DFE (well defined by Assumption~\ref{ass:model}(A2)),
\[
\mathcal R_0=\rho\big(FV^{-1}\big)
=B\sqrt{\frac{\beta_{hm}\,\beta_{mh}\,N_m^\ast}{N\,\mu\,\gamma}}
=B\sqrt{\frac{\beta_{hm}\,\beta_{mh}}{\mu\,\gamma}\cdot m^\ast},
\]
This is the square root of the classical Ross--Macdonald quantity $m^\ast B^2\beta_{hm}\beta_{mh}/(\mu\gamma)$ (without an extrinsic-incubation factor); see, e.g., \cite{RossMacdonaldReview}. The square root is the usual consequence of the next-generation convention, in which a human-to-human cycle through the vector counts as two generations. Both expressions cross the value one simultaneously, so they define the same threshold.
{Here \(m^\ast\) is the mosquito-to-human ratio at the DFE. The linear dependence of \(\mathcal R_0\) on the biting rate \(B\), or equivalently the quadratic dependence of \(\mathcal R_0^2\) on \(B\), together with the proportionality to the vector density \(m^\ast\), is the signature of vector-borne transmission.
}

\medskip
The threshold depends on the disease-free equilibrium at which invasion is assessed. Linearizing at a partially susceptible equilibrium $(S^\ast,0,R^\ast,N_m^\ast,0)$ with $S^\ast\in[0,N]$ replaces the entry $\beta_{hm}B$ of $F$ by $\beta_{hm}B\,S^\ast/N$ (the susceptible human fraction), while the vector entry $\beta_{mh}B\,N_m^\ast/N$ is unchanged; hence the corresponding spectral radius is
\[
\mathcal R_0(S^\ast)=\sqrt{\frac{S^\ast}{N}}\;\mathcal R_0,
\qquad
\mathcal R_0=\mathcal R_0(N)=B\sqrt{\frac{\beta_{hm}\,\beta_{mh}}{\mu\,\gamma}\cdot m^\ast}.
\]
The fully susceptible value $\mathcal R_0$ is therefore the worst-case (maximal) threshold, attained at $S^\ast=N$, and is the quantity used throughout. For a single-host SIR model the analogous dependence on the susceptible level is linear rather than square-root, because a single infection type contributes one factor $S^\ast/N$ instead of the geometric mean over the human--vector cycle.

\begin{proposition}[Transmission threshold at the disease-free equilibrium]
\label{prop:R0}
{Assume \(b>\mu\), \(N>0\), and constant parameters.} At the fully susceptible disease-free equilibrium $E_0$ of \eqref{eq:classic},
the infected subsystem $(I,I_m)$ is linearly stable if $\mathcal R_0<1$ and unstable if $\mathcal R_0>1$. Equivalently,
$\mathcal R_0$ is the invasion threshold: an introduced infection initially decays when $\mathcal R_0<1$ and grows
when $\mathcal R_0>1$.
\end{proposition}

{
	
\begin{proof}
	Let \(z=(I,I_m)^T\) denote the vector of infected variables. Linearizing the
	infected subsystem of \eqref{eq:classic} around the disease-free equilibrium
	\[
	E_0=(N,0,0,N_m^\ast,0),
	\qquad
	N_m^\ast=K\left(1-\frac{\mu}{b}\right),
	\]
	we obtain
	\[
	z'=(F-V)z,
	\]
	where \(F\) contains the terms corresponding to the appearance of new
	infections and \(V\) contains the transition and removal terms. More precisely,
	\[
	F=
	\begin{pmatrix}
		0 & \beta_{hm}B\\[10pt]
		\dfrac{\beta_{mh}B\,N_m^\ast}{N} & 0
	\end{pmatrix},
	\qquad
	V=
	\begin{pmatrix}
		\gamma & 0\\[4pt]
		0 & \mu
	\end{pmatrix}.
	\]
	Thus the next-generation matrix is
	\[
	FV^{-1}
	=
	\begin{pmatrix}
		0 & \dfrac{\beta_{hm}B}{\mu}\\[10pt]
		\dfrac{\beta_{mh}B\,N_m^\ast}{\gamma N} & 0
	\end{pmatrix}.
	\]
	Its eigenvalues are
	\[
	\lambda_{\pm}
	=
	\pm
	B
	\sqrt{
		\frac{\beta_{hm}\beta_{mh}N_m^\ast}{N\mu\gamma}
	}.
	\]
	Therefore
	\[
	\rho(FV^{-1})
	=
	B
	\sqrt{
		\frac{\beta_{hm}\beta_{mh}N_m^\ast}{N\mu\gamma}
	}
	=\mathcal R_0.
	\]
	The stability of the infected subsystem is read directly from the matrix
	\[
	A=F-V=
	\begin{pmatrix}
	-\gamma & \beta_{hm}B\\[8pt]
	\dfrac{\beta_{mh}B\,N_m^\ast}{N} & -\mu
	\end{pmatrix}.
	\]
	Indeed, \(\operatorname{tr}(A)=-(\gamma+\mu)<0\) and
	\[
	\det(A)=\gamma\mu-
	\frac{\beta_{hm}\beta_{mh}B^2N_m^\ast}{N}
	=\gamma\mu\,(1-\mathcal R_0^2).
	\]
	If \(\mathcal R_0<1\), then \(\det(A)>0\) and \(\operatorname{tr}(A)<0\), so both eigenvalues of \(A\) have negative real part. If \(\mathcal R_0>1\), then \(\det(A)<0\), so \(A\) has one positive and one negative real eigenvalue. Hence an introduced infection initially decays if \(\mathcal R_0<1\) and initially grows if \(\mathcal R_0>1\), which proves the threshold assertion.
\end{proof}

}

\begin{remark}
The statement concerns the infected subsystem only. Because the human component has no demography or loss of immunity, $E_0$ belongs to a continuum of equilibria and is not hyperbolic for the full system: the disease-free dynamics have a zero eigenvalue at $E_0$. The hyperbolicity hypothesis (A5) of \cite{vandenDriessche2002} is therefore not satisfied, and the general theorem of that reference is not invoked; the direct computation above replaces it.
\end{remark}

\subsection{Digital Twin Implementation}

Following the initial mosquito--human transmission model in equation \eqref{eq:classic}, we formally define its digital twin by introducing a slow time variable $t_s$ that accounts for gradual changes in epidemiological conditions \cite{AFTN}. The digital twin framework considers two time scales:
\begin{itemize}
\item Fast time $t$: the time variable of the transmission dynamics.
\item Slow time $t_s$: the scale at which data arrive and the parameters are re-estimated, reflecting interventions, climate or vector control. In the numerical experiment the updates are monthly.
\end{itemize}

In accordance with the baseline model \eqref{eq:classic}, the total human
population \(N=S+I+R\) is assumed to be constant along the fast-time
dynamics and to be preserved by the data-assimilation corrections. Thus \(N\)
is treated as a fixed population parameter in the digital twin equations, not
as an additional time-dependent state variable. The mosquito total remains
\(N_m(t,t_s)=S_m(t,t_s)+I_m(t,t_s)\), and the incidence terms keep the same
frequency-dependent denominator \(N\) as in the reference model.

	The digital twin system is given by:
	\begin{equation}\label{eq:malaria-digital-twin}
		\begin{cases}
			\dfrac{\partial S}{\partial t}(t,t_s)
			=
			-\beta_{hm}(t_s)\, B(t_s)\, S(t,t_s)
			\dfrac{I_m(t,t_s)}{N},
			\\[6pt]
			\dfrac{\partial I}{\partial t}(t,t_s)
			=
			\beta_{hm}(t_s)\, B(t_s)\, S(t,t_s)
			\dfrac{I_m(t,t_s)}{N}
			-
			\gamma(t_s)\, I(t,t_s),
			\\[6pt]
			\dfrac{\partial R}{\partial t}(t,t_s)
			=
			\gamma(t_s)\, I(t,t_s),
			\\[6pt]
			\dfrac{\partial S_m}{\partial t}(t,t_s)
			=
			b(t_s)
			\left(
			1-\dfrac{N_m(t,t_s)}{K}
			\right)
			N_m(t,t_s)
			\\[6pt] \qquad \qquad-
			\beta_{mh}(t_s)\, B(t_s)\, S_m(t,t_s)
			\dfrac{I(t,t_s)}{N}
			-
			\mu(t_s)\, S_m(t,t_s),
			\\[6pt]
			\dfrac{\partial I_m}{\partial t}(t,t_s)
			=
			\beta_{mh}(t_s)\, B(t_s)\, S_m(t,t_s)
			\dfrac{I(t,t_s)}{N}
			-
			\mu(t_s)\, I_m(t,t_s).
		\end{cases}
	\end{equation}
The denominator \(N\) is the same as in the baseline model \eqref{eq:classic}; by Assumption~\ref{ass:model}(A2) it is a fixed positive constant, so the incidence terms are well defined without any additive regularization.

\medskip

The system receives monthly data at discrete times determined by $t_s$. We assume that $t_s \in T_s$, where $T_s \subset [0,T]$ is a finite set of slow time instants. Let $T_s = \{t_s^k\}_{k=0}^{J}$ with $t_s^0=0$ and $t_s^k < t_s^{k+1}$ for $k=0,\dots,J-1$. If \(t_s^J<T\), we set \(t_s^{J+1}:=T\) for notational convenience; the parameters fitted at \(t_s^J\) are then used on \((t_s^J,T]\). At each time $t_s^k$, with $k \geq 1$, the system receives data that can be used to estimate the parameters $\beta_{hm}(t_s^k)$, $\beta_{mh}(t_s^k)$, $\gamma(t_s^k)$, and $b(t_s^k)$, and to correct the evolution of the digital twin. In the numerical experiment, $B$ and $\mu$ are kept fixed rather than recalibrated. In the constant-human-population setting considered here, these state corrections are assumed to preserve the total human population, namely
\[
\widehat S^{\,k}+\widehat I^{\,k}+\widehat R^{\,k}=N,
\qquad k=1,\dots,J.
\]

\medskip

For each fixed $t_s^k$, the dynamics in $t$ are governed by constant parameters, which are updated incrementally as $t_s$ progresses. This decoupling assumption allows the model to adapt to changing conditions while maintaining computational tractability. At time $t_s^0$, we assume initial conditions $(S_0,I_0,R_0,S_{m,0},I_{m,0})$ are known, together with initial parameter estimates.

\medskip
The digital twin can be interpreted as a multivalued map \cite{AFTN}:
\[
DT : t \in [0,T] \to DT(t) \subset \bigsqcup_{k \geq 0} \left\{(t_{s}^{k}) \times \mathbb{R}^{5}\right\},
\]
where for $t \in (t_{s}^{j}, t_{s}^{j+1}]$, 
\[
DT(t) = \left\{\left(t_{s}^{k}, (S(t,t_{s}^{k}), I(t,t_{s}^{k}), R(t,t_{s}^{k}), S_{m}(t,t_{s}^{k}), I_{m}(t,t_{s}^{k}))\right), \; k=0,1,\dots,j\right\},
\]
and each quintuple is a solution to system \eqref{eq:classic} on $[t_{s}^{k}, T]$ with parameters fixed at their values at $t_{s}^{k}$.

The \emph{Main Digital Twin} (MDT) is defined as the best forecast at each time $t$, corresponding to the branch associated with the closest previous $t_s^k$:
\[
\mathrm{MDT}(t) = \mathbf{x}(t; t_s^k), \quad t \in (t_s^k, t_s^{k+1}],
\]
where $\mathbf{x}(\cdot; t_s^k)$ is the solution of \eqref{eq:classic} with parameters fitted at $t_s^k$ and initial state corrected to the available monthly data at $t_s^k$. 

\medskip

	The previous construction naturally leads to a system with impulses at the
	slow-time instants \(t_s^k\), which can be rigorously described by means of
	Stieltjes derivatives \cite{AFTN}. Indeed, the Main Digital Twin is not merely
	a piecewise reparametrization of the classical model: at each update time
	\(t_s^k\), the state of the digital model is corrected using the newly received
	monthly data. Therefore, the trajectory of the MDT may have discontinuities at
	the instants \(t_s^k\), and these discontinuities play the role of impulses.
	
	Let us introduce the vector notation
	\[
	\mathbf{x}(t)
	=
	\bigl(S(t),I(t),R(t),S_m(t),I_m(t)\bigr)^{\top}
	\in \mathbb{R}^{5}.
	\]
	For each interval \((t_s^k,t_s^{k+1}]\), the MDT is given by the solution
	computed with the parameter values fitted at the slow time \(t_s^k\). We denote
	these values by
	\[
	\Theta_k
	=
	\bigl(
	\beta_{hm}(t_s^k),
	\beta_{mh}(t_s^k),
	B(t_s^k),
	\gamma(t_s^k),
	b(t_s^k),
	\mu(t_s^k)
	\bigr).
	\]
	Accordingly, we define the piecewise vector field
	\[
	\mathcal{F}(t,\mathbf{x})
	=
	\mathcal{F}_k(\mathbf{x}),
	\qquad
	t\in(t_s^k,t_s^{k+1}],
	\]
	where \(\mathcal{F}_k=(\mathcal{F}_{k,1},\dots,\mathcal{F}_{k,5})\) is the
	right-hand side of the malaria system \eqref{eq:classic} with the parameters
	fixed at the values \(\Theta_k\). Here \(N\) denotes the fixed total human population, while
	\[
	N_m=S_m+I_m.
	\]
	Thus
	\[
	\mathcal{F}_k(\mathbf{x})
	=
	\begin{pmatrix}
		-\beta_{hm}^k B^k S\dfrac{I_m}{N}
		\\[6pt]
		\beta_{hm}^k B^k S\dfrac{I_m}{N}-\gamma^k I
		\\[6pt]
		\gamma^k I
		\\[6pt]
		b^k\left(1-\dfrac{N_m}{K}\right)N_m
		-\beta_{mh}^k B^k S_m\dfrac{I}{N}
		-\mu^k S_m
		\\[6pt]
		\beta_{mh}^k B^k S_m\dfrac{I}{N}
		-\mu^k I_m
	\end{pmatrix}.
	\]
	Here, for instance, \(\beta_{hm}^k=\beta_{hm}(t_s^k)\),
	\(\beta_{mh}^k=\beta_{mh}(t_s^k)\), and analogously for the remaining
	parameters.
	
	At the update time \(t_s^k\), \(k\geq1\), the available monthly data determine a
	\emph{correction map} \(\mathcal C_k:\mathbb{R}^5\to\mathbb{R}^5\), which sends the
	pre-update state \(\mathbf{x}\) to the corrected state
	\[
	\widehat{\mathbf{x}}^k
	=
	\mathcal C_k(\mathbf{x})
	=
	\bigl(
	\widehat S^k,
	\widehat I^k,
	\widehat R^k,
	\widehat S_m^k,
	\widehat I_m^k
	\bigr)^{\top}.
	\]
	Two cases are relevant. If every compartment is observed, \(\mathcal C_k\) is
	constant (a full reset to the data). In the numerical experiment only \(I\) and
	the mosquito compartments are tied to the observations, whereas the recovered
	compartment is carried from the model:
	\begin{equation}\label{eq:correction-map}
	\widehat I^k=\frac{Y_k}{q},
	\qquad
	\widehat R^k=\min\bigl\{\max\{R,0\},\,N-\widehat I^k\bigr\},
	\qquad
	\widehat S^k=N-\widehat I^k-\widehat R^k,
	\end{equation}
	with \((\widehat S_m^k,\widehat I_m^k)\) prescribed by the mosquito observations.
	In both cases \(\mathcal C_k\) is globally Lipschitz (the truncation in
	\eqref{eq:correction-map} is \(1\)-Lipschitz) and, consistently with the
	constant-human-population assumption,
	\[
	\widehat S^k+\widehat I^k+\widehat R^k=N,
	\qquad k=1,\dots,J.
	\]
	If, in addition, \(0\le Y_k/q\le N\), then \(\widehat S^k,\widehat I^k,\widehat R^k\ge0\).
	Throughout, we only assume that each \(\mathcal C_k\) is globally Lipschitz and
	takes values in the hyperplane \(\{S+I+R=N\}\). The impulse at \(t_s^k\) is therefore
	\[
	\mathcal{J}_k(\mathbf{x})
	=
	\mathcal C_k(\mathbf{x})-\mathbf{x}.
	\]
	Consequently, the MDT can be written as the solution of the impulsive system
	\[
	\begin{cases}
		\mathbf{x}'(t)=\mathcal{F}(t,\mathbf{x}(t)),
		& \text{for a.e. } t\in[0,T]\setminus \{t_s^1,\dots,t_s^J\},
		\\[4pt]
		\mathbf{x}(0)=\mathbf{x}_0,
		\\[4pt]
		\mathbf{x}\bigl((t_s^k)^+\bigr)
		=
		\mathbf{x}(t_s^k)+\mathcal{J}_k(\mathbf{x}(t_s^k))
		=
		\mathcal C_k(\mathbf{x}(t_s^k)),
		& k=1,\dots,J.
	\end{cases}
	\tag{MDT-I}
	\label{eq:MDT-impulsive}
	\]
	Thus, the MDT follows the classical malaria dynamics between two consecutive
	data updates, whereas at the update times it is instantaneously corrected
	towards the state inferred from the available data.
	
	We now rewrite \eqref{eq:MDT-impulsive} as a Stieltjes differential equation.
	Let \(\alpha_k>0\), \(k=1,\dots,J\), and define the derivator
	\[
	g_{DT}(t)
	=
	t+\sum_{\substack{1\le k\le J\\ t_s^k<t}}\alpha_k,
	\qquad t\in[0,T].
	\]
	The function \(g_{DT}\) is nondecreasing and left-continuous. Moreover,
	\[
	\Delta g_{DT}(t_s^k)
	=
	g_{DT}\bigl((t_s^k)^+\bigr)-g_{DT}(t_s^k)
	=
	\alpha_k,
	\]
	so that the associated Lebesgue--Stieltjes measure is
	\[
	\mu_{g_{DT}}
	=
	\mu_{Id}
	+
	\sum_{k=1}^{J}\alpha_k\delta_{t_s^k}.
	\]
	With the left-continuous convention used in this paper, \(\mathbf{x}(t_s^k)\) denotes the pre-update state and \(\mathbf{x}((t_s^k)^+)\) denotes the state immediately after correction. Integrals over \([0,t)\) include precisely the atoms whose locations are strictly smaller than \(t\).
	The continuous part of this measure accounts for the usual chronological
	evolution of the malaria model, while the atomic part encodes the corrections
	of the digital twin at the data assimilation times.
	
	Define
	\[
	\mathcal{F}_{DT}(t,\mathbf{x})
	=
	\begin{cases}
		\mathcal{F}(t,\mathbf{x}),
		& t\notin \{t_s^1,\dots,t_s^J\},
		\\[6pt]
		\displaystyle
		\frac{\mathcal{J}_k(\mathbf{x})}{\alpha_k},
		& t=t_s^k,\quad k=1,\dots,J.
	\end{cases}
	\]
	Then the impulsive problem \eqref{eq:MDT-impulsive} is equivalent to the
	Stieltjes differential system
	\[
	\begin{cases}
		\mathbf{x}'_{g_{DT}}(t)
		=
		\mathcal{F}_{DT}(t,\mathbf{x}(t)),
		& \mu_{g_{DT}}\text{-a.e. } t\in[0,T],
		\\[4pt]
		\mathbf{x}(0)=\mathbf{x}_0.
	\end{cases}
	\tag{MDT-\(g\)}
	\label{eq:MDT-Stieltjes}
	\]
	Indeed, outside the update times the derivator satisfies \(d\mu_{g_{DT}}=d\mu_{Id}=dt\), and
	therefore the \(g_{DT}\)-derivative coincides with the classical derivative.
	At an update time \(t_s^k\), however, the Stieltjes derivative gives
	\[
	\mathbf{x}'_{g_{DT}}(t_s^k)
	=
	\frac{
		\mathbf{x}\bigl((t_s^k)^+\bigr)-\mathbf{x}(t_s^k)
	}{
		g_{DT}\bigl((t_s^k)^+\bigr)-g_{DT}(t_s^k)
	}
	=
	\frac{\mathcal{J}_k(\mathbf{x}(t_s^k))}{\alpha_k}.
	\]
	Hence,
	\[
	\mathbf{x}\bigl((t_s^k)^+\bigr)
	=
	\mathbf{x}(t_s^k)
	+
	\alpha_k\, \mathcal{F}_{DT}(t_s^k,\mathbf{x}(t_s^k))
	=
	\mathcal C_k(\mathbf{x}(t_s^k)),
	\]
	which is precisely the correction step defining the MDT.
	
	This formulation gives a precise meaning to the Main Digital Twin as a
	\(g_{DT}\)-absolutely continuous function
	\[
	\mathbf{x}\in \mathcal{AC}_{g_{DT}}\bigl([0,T];\mathbb{R}^5\bigr)
	\]
	satisfying \eqref{eq:MDT-Stieltjes}. In integral form, this reads
	\[
	\mathbf{x}(t)
	=
	\mathbf{x}_0
	+
	\int_{[0,t)}
	\mathcal{F}_{DT}(s,\mathbf{x}(s))\,d\mu_{g_{DT}}(s).
	\]
	Equivalently,
	\[
	\mathbf{x}(t)
	=
	\mathbf{x}_0
	+
	\int_0^t
	\mathcal{F}(s,\mathbf{x}(s))\,ds
	+
	\sum_{\substack{1\le k\le J\\ t_s^k<t}}
	\mathcal{J}_k(\mathbf{x}(t_s^k)).
	\]
	Thus, the Stieltjes formulation unifies in a single equation the continuous
	malaria transmission dynamics and the discrete data-assimilation corrections
	that define the digital twin.
	
	From the analytical point of view, this interpretation is useful because the
	existence of the MDT can be studied using the theory of systems of Stieltjes
	differential equations. Since \(T_s=\{t_s^0,\dots,t_s^J\}\) is finite and the
	parameters are piecewise constant, the derivator \(g_{DT}\) has only finitely
	many jumps. Moreover, since \(N>0\) is a fixed constant, the vector field \(\mathcal{F}\)
	is polynomial, hence locally Lipschitz, in the epidemiological variables, and the correction maps
	\(\mathcal{J}_k\) are globally Lipschitz. Therefore, standard local existence and uniqueness
	results for Stieltjes differential systems apply. Consequently, there exists
	\(\tau\in(0,T]\) such that the Stieltjes problem
	\eqref{eq:MDT-Stieltjes} has a unique solution on \([0,\tau]\). If, in addition,
	the solution remains bounded, it can be continued up to the final time \(T\).
	
	\begin{theorem}[Local well-posedness of the Main Digital Twin]
\label{thm:MDT-existence-uniqueness}
			Let
			\[
			T_s=\{t_s^0,\dots,t_s^J\},
			\qquad
			0=t_s^0<t_s^1<\cdots<t_s^J\leq T,
			\]
			be the finite set of update times of the Main Digital Twin. If \(t_s^J<T\), set \(t_s^{J+1}:=T\) and use the last fitted parameter values on \((t_s^J,T]\). The continuous vector field is therefore understood on the intervals \((t_s^k,t_s^{k+1}]\) for all relevant \(k\). Let
			\[
			g_{DT}(t)
			=
			t+\sum_{\substack{1\le j\le J\\ t_s^j<t}}\alpha_j,
			\qquad
			\alpha_j>0,
			\]
			be the associated digital-twin derivator. Consider the Stieltjes differential
			problem
			\[
			\begin{cases}
				\mathbf{x}'_{g_{DT}}(t)
				=
				\mathcal{F}_{DT}(t,\mathbf{x}(t)),
				& \mu_{g_{DT}}\text{-a.e. } t\in[0,T],
				\\[4pt]
				\mathbf{x}(0)=\mathbf{x}_0,
			\end{cases}
			\]
			where
			\[
			\mathbf{x}(t)
			=
			(S(t),I(t),R(t),S_m(t),I_m(t))^\top.
			\]
			The constant-human-population formulation is first considered on the affine hyperplane
			\[
			\mathcal A_N=\{\mathbf{x}\in\mathbb R^5:\ S+I+R=N\}.
			\]
			Positivity is a separate invariance property and is established under the additional biological hypotheses of Corollary~\ref{cor:MDT-global}.
			Assume that the parameter functions are piecewise constant with respect to the
			partition induced by \(T_s\), and that the impulses at the update times
			are given by
			\[
			\mathcal{J}_j(\mathbf{x})
			=
			\mathcal C_j(\mathbf{x})-\mathbf{x},
			\qquad j=1,\dots,J,
			\]
			where each correction map \(\mathcal C_j:\mathbb R^5\to\mathbb R^5\) is globally Lipschitz and preserves the fixed human population, namely
			\[
			\mathcal C_j(\mathcal A_N)\subset\mathcal A_N,
			\qquad j=1,\dots,J.
			\]
			Then, for every initial condition
			\[
			\mathbf{x}_0\in\mathbb{R}_+^5,
			\qquad
			S_0+I_0+R_0=N,
			\]
			there exists \(\tau\in(0,T]\) such that problem
			\eqref{eq:MDT-Stieltjes} has a unique solution
			\[
			\mathbf{x}\in
			\mathcal{AC}_{g_{DT}}\bigl([0,\tau];\mathbb{R}^5\bigr).
			\]
			Under the population-preserving assumptions above, this solution satisfies \(\mathbf{x}(t)\in\mathcal A_N\) for all \(t\in[0,\tau]\).
			Moreover, if the solution remains bounded, then it can be continued up to the final time \(T\).
		\end{theorem}
		
		\begin{proof}
			We verify that the problem satisfies the usual local existence and uniqueness
			hypotheses for first-order systems of Stieltjes differential equations, as
			developed in \cite{FrigonLopezPouso2017}. The interpretation of impulsive
			differential equations as Stieltjes differential equations is also used in
			\cite{LopezPousoRodriguez2015,AFTN}, and general existence principles with
			applications to biological models can be found in
			\cite{LopezPousoMarquez2018}.
			
			First, \(g_{DT}\) is a derivator. Indeed, since each \(\alpha_j>0\), the function
			\(g_{DT}\) is nondecreasing and left-continuous. Its associated
			Lebesgue--Stieltjes measure is
			\[
			\mu_{g_{DT}}
			=
			\mu_{Id}
			+
			\sum_{j=1}^{J}\alpha_j\delta_{t_s^j}.
			\]
			Thus, \(\mu_{g_{DT}}\) is a finite positive measure on \([0,T]\), with an
			absolutely continuous part with respect to the classical Lebesgue measure and a
			finite atomic part concentrated at the update times.
			
			Second, for every fixed state \(\mathbf{x}\in\mathbb{R}^5\), the map
			\[
			t\mapsto \mathcal{F}_{DT}(t,\mathbf{x})
			\]
			is \(\mu_{g_{DT}}\)-measurable. Outside the finite set \(T_s\), the vector field
			coincides with the malaria vector field whose parameters are frozen on each
			interval \((t_s^j,t_s^{j+1}]\). At an update time \(t_s^j\), it is prescribed by
			\[
			\mathcal{F}_{DT}(t_s^j,\mathbf{x})
			=
			\frac{\mathcal{J}_j(\mathbf{x})}{\alpha_j}.
			\]
			Therefore, the time dependence is piecewise constant, up to finitely many
			distinguished points, and hence it is \(\mu_{g_{DT}}\)-measurable.
			
			Third, the map
			\[
			\mathbf{x}\mapsto \mathcal{F}_{DT}(t,\mathbf{x})
			\]
			is Lipschitz on bounded sets. For \(M>0\) let
			\[
			B_M
			=
			\left\{
			\mathbf{x}\in\mathbb{R}^5:\ \|\mathbf{x}\|\le M
			\right\}.
			\]
			No sign condition is imposed at this stage: balls are used so that every initial condition, including those on the boundary of \(\mathbb{R}_+^5\), has a neighbourhood on which the estimates hold; nonnegativity is recovered under the hypotheses of Corollary~\ref{cor:MDT-global}.
			Since \(N\) is a fixed positive constant (Assumption~\ref{ass:model}(A2)), the continuous part of the malaria vector field is polynomial in \(\mathbf{x}\), with bounded piecewise constant coefficients, and is therefore Lipschitz on \(B_M\). At the update times, the
			maps \(\mathcal{J}_j=\mathcal C_j-\mathrm{Id}\) are globally Lipschitz, because the correction maps \(\mathcal C_j\) are. Since the number of update times
			is finite, there exists a constant \(L_M>0\) such that
			\[
			\|\mathcal{F}_{DT}(t,\mathbf{x})
			-\mathcal{F}_{DT}(t,\mathbf{y})\|
			\leq
			L_M\|\mathbf{x}-\mathbf{y}\|,
			\]
			for every \(\mathbf{x},\mathbf{y}\in B_M\) and for every \(t\in[0,T]\).
			
			Fourth, the right-hand side is locally bounded by a function in
			\(L^1_{g_{DT}}([0,T])\). Indeed, on the continuous part this follows from the
			boundedness of the parameters on each interval \((t_s^j,t_s^{j+1}]\) and from
			the boundedness of \(\mathbf{x}\in B_M\). At the atoms, it follows from
			\[
			\mathcal{F}_{DT}(t_s^j,\mathbf{x})
			=
			\frac{\mathcal C_j(\mathbf{x})-\mathbf{x}}{\alpha_j},
			\]
			from the continuity of \(\mathcal C_j\), and from the finiteness of \(T_s\). Hence, for every \(M>0\) there exists a constant \(h_M\) (in particular \(h_M\in L^1_{g_{DT}}([0,T])\), the measure being finite) such that
			\[
			\|\mathcal{F}_{DT}(t,\mathbf{x})\|
			\leq h_M,
			\qquad
			\mathbf{x}\in B_M,\quad t\in[0,T].
			\]
			
			Therefore, the hypotheses of the local existence and uniqueness theorem for
			Stieltjes differential systems are satisfied. Consequently, there exists
			\(\tau\in(0,T]\) and a unique solution
			\[
			\mathbf{x}\in
			\mathcal{AC}_{g_{DT}}\bigl([0,\tau];\mathbb{R}^5\bigr)
			\]
			of \eqref{eq:MDT-Stieltjes}. Equivalently, this solution satisfies the
			integral equation
			\[
			\mathbf{x}(t)
			=
			\mathbf{x}_0
			+
			\int_{[0,t)}
			\mathcal{F}_{DT}(s,\mathbf{x}(s))\,d\mu_{g_{DT}}(s),
			\qquad t\in[0,\tau].
			\]
			
			The solution remains in \(\mathcal A_N\): between update times the sum of the three human equations vanishes, and at an update time \(\mathbf{x}((t_s^j)^+)=\mathcal C_j(\mathbf{x}(t_s^j))\in\mathcal A_N\) by hypothesis.

			Finally, if the solution remains bounded, say in \(B_M\),
			then the same argument can be restarted from any point of the interval of
			existence, with constants depending only on \(M\). Since \([0,T]\) is bounded and \(T_s\) is finite, the solution can be
			continued step by step up to \(T\). This proves the claimed continuation
			criterion.
		\end{proof}

\begin{corollary}[Global well-posedness of the Main Digital Twin]
\label{cor:MDT-global}
Let
\[
\Omega^\ast
=
\Big\{\mathbf{x}\in\mathbb{R}_+^5:
S+I+R=N,\quad S_m+I_m\le K
\Big\}.
\]
Suppose, in addition to the hypotheses of Theorem \ref{thm:MDT-existence-uniqueness}, that the
parameters satisfy the standing assumptions of Lemma \ref{lem:positivity} on every inter-update
interval (in particular, \(b^k>\mu^k\)), that the initial datum is biologically admissible,
\(\mathbf{x}_0\in\Omega^\ast\), and that the correction maps preserve the admissible region,
\[
\mathcal C_k(\Omega^\ast)\subset\Omega^\ast,
\qquad k=1,\dots,J.
\]
Then the unique local solution of \eqref{eq:MDT-Stieltjes} extends to the whole
interval \([0,T]\) and remains in \(\Omega^\ast\); that is, the Main Digital Twin is a well-defined \(g_{DT}\)-absolutely continuous
function on \([0,T]\).
\end{corollary}

For a full reset, the condition \(\mathcal C_k(\Omega^\ast)\subset\Omega^\ast\) means that the assimilated state \(\widehat{\mathbf{x}}^{\,k}\) is nonnegative, has human total \(N\) and satisfies \(\widehat S_m^{\,k}+\widehat I_m^{\,k}\le K\). For the correction \eqref{eq:correction-map} it holds whenever \(0\le Y_k/q\le N\) and the mosquito observations are nonnegative with total at most \(K\).

\begin{proof}
By Lemma \ref{lem:positivity}, applied on each interval \((t_s^{k},t_s^{k+1}]\), where the parameters
are constant and satisfy the standing assumptions, a trajectory that starts in \(\Omega^\ast\) remains in \(\Omega^\ast\) between two
consecutive update times. At each update time the state is mapped to \(\mathcal C_k(\mathbf{x}(t_s^k))\), which belongs to \(\Omega^\ast\) by hypothesis. Hence the solution stays
in the compact set \(\Omega^\ast\) on its whole interval of existence. In particular it remains bounded, and the continuation criterion of Theorem \ref{thm:MDT-existence-uniqueness} gives existence up to \(t=T\).
\end{proof}

In this sense, the MDT associated with the malaria transmission model is not only a numerical forecasting procedure, but also the solution of a well-defined Stieltjes differential equation. The atoms of the derivator encode the moments at which the digital twin interacts with the observed system, while the values of \(\mathcal{F}_{DT}\) at those atoms encode the magnitude and direction of the corrections induced by the incoming data. This is especially appropriate for mosquito-borne diseases, where the vector component may be strongly affected by intermittent environmental information, surveillance data or sudden control interventions.

\subsection{Digital Twin Implementation in Python}

The numerical experiments were carried out with a self-contained Python implementation based on \texttt{numpy}, \texttt{scipy}, \texttt{pandas} and \texttt{matplotlib}. The input contains 132 monthly observations of confirmed uncomplicated malaria in Kebbi State, from January 2015 to December 2025; the window January 2015--May 2024 is the period covered by the retrospective Kebbi case series analysed in Bakare et al.~\cite{Bakare2026}, and monthly counts are used throughout. The monthly records are internally consistent: confirmed uncomplicated cases agree with the sum of positive rapid-diagnostic and microscopy tests (median ratio $0.98$). Two months are anomalously low, January 2015 ($177$ cases) and May 2018 ($9{,}143$); in both, the number of persons tested falls in the same proportion, which is consistent with incomplete reporting but does not establish its cause; both observations are kept unmodified. The source also contains annual aggregates; they coincide with the sums of the monthly values except for 2023, where the annual totals of the case, testing and treatment indicators exceed the monthly sums by $10$--$17\%$ (for confirmed cases, $921{,}714$ against $817{,}576$), while the total of persons with fever coincides. Only the monthly values are used here. No observed entomological series accompanies these data; therefore, the mosquito series used to compare the classical and Stieltjes formulations is synthetic and is generated as described in Section~\ref{sec:secondmodel}. All numerical times and rates are expressed in months. The fixed demographic anchors are $N=4{,}377{,}770$, $K=8.0\times10^7$ and $\mu=2.14\,\mathrm{month}^{-1}$, consistently with the data context of Bakare et al.~\cite{Bakare2026}.

The reported monthly cases are incidence counts, whereas the state variable $I(t)$ is an instantaneous infectious compartment. To make this distinction explicit, the numerical experiment uses the preliminary observation equation
\[
Y_k=qI(t_s^k)+\varepsilon_k,
\]
with configurable scale $q>0$ and default value $q=1$. Thus the case series is treated as a proxy for $I$, not as a validated prevalence reconstruction. A fully epidemiological calibration would require estimating $q$, introducing reporting and delay mechanisms, or fitting integrated monthly incidence.

The implementation uses the following settings:
\begin{itemize}
\item the first seven months form the baseline period and the common evaluation begins after it;
\item parameters are estimated on rolling four-month windows by bounded nonlinear least squares with Tikhonov regularization toward the preceding estimate;
\item only $(\beta_{hm},\beta_{mh},\gamma,b)$ are fitted, while $B$, $\mu$, $K$ and the seasonal/impulse parameters are fixed;
\item each accepted estimate is smoothed as $0.7\,\Theta_{\rm new}+0.3\,\Theta_{\rm old}$;
\item at month $t_s^k$ the state is corrected by the map \eqref{eq:correction-map}: $I$ is set to $Y_k/q$, $R$ is carried from the model and truncated to $[0,N-Y_k/q]$, and $S=N-I-R$; the mosquito total and infected fraction are set to the synthetic observations;
\item both the classical and Stieltjes twins are run by the same driver, so their errors are calculated on an identical evaluation set;
\item continuous segments are integrated by fixed-step RK4. In the Stieltjes model the biological atoms are integer month indices, so they coincide with update times (convention (II) of Section~\ref{sec:secondmodel}): the impulse is applied when the continuous evolution reaches $\tau_k$, the resulting pre-assimilation state is recorded as the model value at that date, and the observation correction is applied afterwards.
\end{itemize}

For both runs, the code also computes the online transmission indicator
\begin{equation}\label{eq:RDT}
\mathcal R_{\mathrm{DT}}(t_s^k)=B^k\sqrt{\frac{\beta_{hm}^k\beta_{mh}^k\,\widehat N_m^{\,k}}{N\mu^k\gamma^k}}.
\end{equation}
This is the next-generation expression of Proposition~\ref{prop:R0} with the disease-free mosquito equilibrium $N_m^*$ replaced by the current corrected mosquito abundance $\widehat N_m^{\,k}$. Consequently, $\mathcal R_{\mathrm{DT}}(t)$ is an adaptive indicator of instantaneous transmission potential; it is not the basic reproduction number of the nonautonomous, repeatedly corrected digital twin. The autonomous-model $\mathcal R_0$ remains the quantity in Proposition~\ref{prop:R0}.

The code exports trajectories, one-step forecasts, uncertainty bands, parameter histories, $\mathcal R_{\mathrm{DT}}$, metrics, the synthetic mosquito data, an Excel workbook and all figures used below. RMSE, MAE and $R^2$ are computed using their standard definitions on
finite, common observations after the baseline.

\subsection{Incremental Fitting, Correction and Forecasting}

At each monthly update time $t_s^k$ from the eighth observation onward, the algorithm performs the following operations:
\begin{enumerate}
\item it propagates the preceding corrected state to $t_s^k$; in the Stieltjes model any biological atom at $t_s^k$ is applied at the end of this propagation;
\item it applies the observation correction at $t_s^k$;
\item it estimates $(\beta_{hm},\beta_{mh},\gamma,b)$ on the preceding four-month window;
\item it starts from the corrected state and issues a genuine one-month-ahead forecast.
\end{enumerate}
The human operational forecast is a convex combination of the median mechanistic ensemble (weight $0.2$) and a past-only seasonal-naive prediction (weight $0.8$). The latter uses only information available at the forecast origin, namely $Y_{k+1-12}$ when a twelve-month lag exists and the most recent observation otherwise. The $10$--$90\%$ ensemble envelope shown in the figures is obtained by perturbing the fitted mechanistic parameters by $15\%$ and then applying the same convex combination; it is not a statistically calibrated predictive interval.

\begin{figure}[ht!]
\centering
\begin{tikzpicture}[
  >=Latex,
  line cap=round,
  nodeMinW/.style={minimum width=5.2cm, minimum height=1.4cm},
  nodeTextW/.style={text width=4.9cm},
  every node/.style={
    draw,
    rectangle,
    rounded corners,
    align=center,
    font=\small
  }
]
% Left column (data flow)
\node[fill=blue!10, nodeMinW, nodeTextW] (data) at (0.0, 0.0) 
    {Data Stream\\(Monthly human cases, synthetic mosquito series)};
    
\node[fill=blue!10, nodeMinW, nodeTextW] (window) at (0.0, -3.2) 
    {Sliding Window Selection\\(4-month fit + forecast)};
    
\node[fill=yellow!20, nodeMinW, nodeTextW] (model) at (0.0, -6.4) 
    {Transmission Model\\Classical or Stieltjes equations};

% Right column (digital twin process)  
\node[fill=orange!20, nodeMinW, nodeTextW] (optim) at (8.0, -1.6) 
    {Parameter Optimization\\(\texttt{scipy}, regularized least squares)};
    
\node[fill=purple!10, nodeMinW, nodeTextW] (dtwin) at (8.0, -4.8) 
    {Digital Twin\\Updated parameters \& corrected state};
    
\node[fill=red!10, nodeMinW, nodeTextW] (pred) at (8.0, -8.0) 
    {Predictions \& Hospital Planning\\Infection waves, bed requirements};

% Bottom centre (export)
\node[fill=gray!10, minimum width=4.5cm, minimum height=1.2cm, draw, rounded corners] (export) at (4.0, -11.0) 
    {Export to Excel\\Parameter history, performance metrics};

% Data flow - vertical arrows
\draw[->, very thick] (data) -- (window) node[midway, left] {Data input};
\draw[->, very thick] (window) -- (model) node[midway, left] {Initialize state};

% Model to optimization
\draw[->, very thick] (model.east) -- ++(1.0,0) |- (optim.west) 
    node[midway, above, yshift=0.8cm, xshift=1.5cm] {Model predictions};

% Optimization and feedback - vertical arrows on right
\draw[->, very thick] (optim.south) -- (dtwin.north) 
    node[midway, right] {Update parameters};
\draw[->, very thick] (dtwin.south) -- (pred.north) 
    node[midway, right] {Generate forecasts};
    
% Parameter feedback loop
\draw[->, very thick, dashed] (dtwin.west) .. controls ($ (dtwin.west) + (-2.5, 0.5) $) and ($ (model.east) + (-0.5, 0.5) $) .. 
    (model.east) node[midway, above, yshift=0.8cm, xshift=-1.0cm] {Parameter feedback};

% New data trigger update
\draw[->, thick, dashed, blue] (data.east) .. controls ($ (data.east) + (2.5, 0.7) $) and ($ (window.east) + (2.5, -0.3) $) .. 
    (window.east) node[midway, above, xshift=1.2cm, yshift=1.4cm] {New data at $t_s^k$};

% Export of results
\draw[->, thick] (pred.south) -- ++(0, -1.0) -| (export.north);

% Error metrics from model to export
\draw[->, thick, dashed] (model.south) .. controls ($ (model.south) + (0, -1.2) $) and ($ (export.west) + (-0.7, 0.7) $) .. 
    (export.west) node[midway, left, xshift=-0.4cm, yshift=0.3cm] {Error metrics};

\end{tikzpicture}
\caption{Digital twin workflow for malaria prediction. Data enter the system at slow time points $t_s^k$, triggering parameter recalibration within 4-month sliding windows. The transmission model (classical or Stieltjes-based) is calibrated by bounded least squares with regularization, updating the digital twin's parameters and correcting states through jump updates. One-month-ahead predictions are generated, together with an illustrative conversion into bed demand, and the results are exported for analysis. The loop is repeated whenever new data arrive.}
\label{fig:workflow-feedback}
\end{figure}
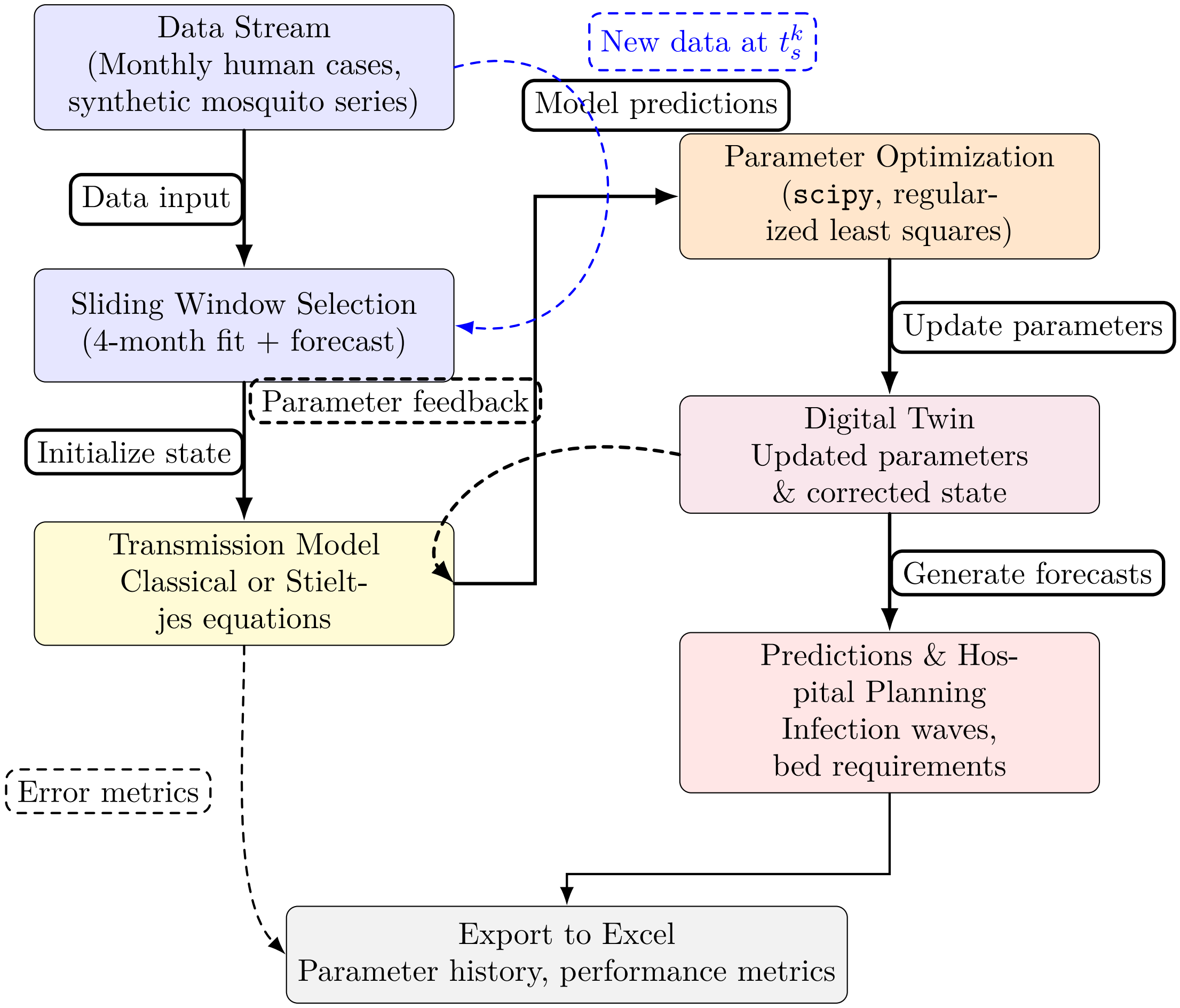

Figure~\ref{fig:workflow-feedback} summarizes the workflow.

\medskip

The classical digital twin of this section supplies the reference implementation. Its mosquito subsystem evolves continuously and therefore cannot represent concentrated biological events such as emergence pulses. This motivates the extension of the digital twin framework by means of Stieltjes differential calculus: replacing the time derivative in the mosquito compartments by a derivative with respect to an appropriate derivator incorporates continuous development and concentrated events in a single model, while preserving the structure of the digital twin and its data-assimilation mechanism. The numerical results for both twins are reported together in Section~\ref{sec:results}.

\section{Digital Twin Model with Stieltjes-Driven Mosquito Dynamics}
\label{sec:secondmodel}

\subsection{Biological Motivation for Stieltjes Derivatives}

The life cycle of \emph{Anopheles} mosquitoes is characterized by impulsive biological events---egg hatching, pupal emergence, and seasonal diapause---that introduce discontinuities in population dynamics. Classical continuous-time derivatives cannot adequately represent these abrupt changes. To address this limitation, we incorporate Stieltjes derivatives into the mosquito compartments of the digital twin model, extending the methodology previously applied to invasive species such as \emph{Vespa velutina} \cite{Vespa} and leveraging digital twin theory based on Stieltjes calculus \cite{AFTN}.

\medskip

The genus \textit{Anopheles} (Diptera: Culicidae) comprises over 400 species, with approximately 30 serving as primary vectors of human malaria \cite{Beier1998,Clements1992}. These mosquitoes undergo holometabolous development through four distinct stages: egg, larva, pupa, and adult. The progression through these stages is highly sensitive to environmental conditions such as temperature, humidity, and breeding site availability \cite{Clements1992}.

\medskip

Key impulsive transitions relevant for epidemiological modelling include:
\begin{itemize}
\item Egg hatching: sudden influxes of larvae following rainfall or temperature thresholds.
\item Pupal emergence: discrete increases in adult mosquito numbers as cohorts complete metamorphosis.
\item Seasonal diapause: rapid population decline during unfavourable conditions and reactivation when conditions improve.
\end{itemize}

These discontinuous biological processes motivate the use of Stieltjes derivatives, which provide a unified framework for combining smooth growth phases with discrete jumps. The Stieltjes derivator $g_m(t)$ encodes the timing and magnitude of these impulses, allowing the digital twin to reproduce realistic mosquito population dynamics.

\subsection{Mathematical Formulation}

In the hybrid formulation, the human compartments evolve with respect to chronological time $g_H(t)=t$, whereas the mosquito compartments evolve with respect to a Stieltjes derivator $g_m$ whose atoms are the biological impulse times $\tau_k$; see \eqref{eq:derivator}. The numerical experiment uses
\[
g_m(t)=t+\sum_{\tau_k<t}\delta_k,
\qquad \delta_k=0.5,
\]
so that the continuous part of $g_m$ is chronological time and all continuous rates retain monthly units. The atoms $\tau_k$ occur once per year, at the rainy-season onset months; in the numerical experiment they are the pulse months of the synthetic generator, which are determined from the seasonal climatology of the complete case series (Section~\ref{sec:numerics}). The jump size $\delta_k$ belongs to the derivator, while the actual population change is prescribed by the impulse map $J_k$ through the atomic field $J_k/\delta_k$; thus changing $\delta_k$ alone does not change the prescribed jump. The fixed human total is denoted by $N$, and the Stieltjes mosquito subsystem includes both annual seasonal recruitment and the discrete emergence pulses.

\subsubsection{Hybrid System with Stieltjes Mosquito Dynamics}

The classical system \eqref{eq:classic} is enhanced by incorporating Stieltjes derivatives in the mosquito dynamics:
\begin{equation}\label{eq:enhanced}
\begin{cases}
\dfrac{dS}{dt}
=
-\beta_{hm} B S \dfrac{I_m}{N}, \\[8pt]
\dfrac{dI}{dt}
=
\beta_{hm} B S \dfrac{I_m}{N}
-
\gamma I, \\[8pt]
\dfrac{dR}{dt}
=
\gamma I, \\[8pt]
\dfrac{dS_m}{dg_m}
=
\mathcal F^{S_m}(t,I,S_m,I_m), \\[8pt]
\dfrac{dI_m}{dg_m}
=
\mathcal F^{I_m}(t,I,S_m,I_m),
\end{cases}
\end{equation}
where \(N_m=S_m+I_m\), \(0\le s_A\le1\) (so that the seasonal recruitment factor is nonnegative), and the Stieltjes vector fields of the mosquito compartments are defined \emph{piecewise},
so that the continuous dynamics act on the absolutely continuous part of $g_m$ and the biological
impulses act only at the atoms $\tau_k$:
\begin{equation}\label{eq:mosquito-field-S}
\mathcal F^{S_m}(t,I,S_m,I_m)
=
\begin{cases}
b\left(1+s_A\sin\left(\dfrac{2\pi t}{T_{\mathrm{season}}}+\phi\right)\right)
\left(1-\dfrac{N_m}{K}\right)N_m\\[10pt]
-\beta_{mh} B S_m \dfrac{I}{N}
-\mu S_m,
& t\notin\{\tau_k\}, \\[10pt]
\dfrac{J^S_k(S_m,I_m)}{\delta_k},
& t=\tau_k,
\end{cases}
\end{equation}
\begin{equation}\label{eq:mosquito-field-I}
\mathcal F^{I_m}(t,I,S_m,I_m)
=
\begin{cases}
\beta_{mh} B S_m \dfrac{I}{N}
-\mu I_m,
& t\notin\{\tau_k\}, \\[10pt]
\dfrac{J^I_k(S_m,I_m)}{\delta_k},
& t=\tau_k.
\end{cases}
\end{equation}

The Stieltjes derivator is
\begin{equation}\label{eq:derivator}
g_m(t)
=
\alpha t
+
\sum_{\{k\in\mathcal K_m:\,\tau_k<t\}}\delta_k,
\end{equation}
with $\alpha>0$, $\delta_k>0$ and a finite index set $\mathcal K_m$. For $\alpha\neq1$ the continuous rates are rescaled by the factor $\alpha$; throughout the analysis and in the numerical experiment we take $\alpha=1$.

The impulse maps are defined by
\begin{equation}\label{eq:mosquito-impulses}
J^S_k(S_m,I_m)
=
\rho^S_k S_m+\sigma_k,
\qquad
J^I_k(S_m,I_m)
=
\rho^I_k I_m.
\end{equation}

In order to preserve the nonnegativity of the mosquito compartments after impulses, we assume
\[
\rho_k^S\ge -1,
\qquad
\rho_k^I\ge -1,
\qquad
\sigma_k\ge0,
\qquad k\in\mathcal K_m.
\]
Indeed, under these conditions,
\[
S_m(\tau_k^+)=(1+\rho_k^S)S_m(\tau_k)+\sigma_k\ge0,
\qquad
I_m(\tau_k^+)=(1+\rho_k^I)I_m(\tau_k)\ge0,
\]
whenever \(S_m(\tau_k),I_m(\tau_k)\ge0\).

With the piecewise definition above, the atomic part of $g_m$ produces, at each impulse time $\tau_k$, the jump
\[
S_m(\tau_k^+)-S_m(\tau_k)
=
\delta_k\,\mathcal F^{S_m}(\tau_k,I(\tau_k),S_m(\tau_k),I_m(\tau_k))
=
J^S_k(S_m(\tau_k),I_m(\tau_k)),
\]
and analogously $I_m(\tau_k^+)-I_m(\tau_k)=J^I_k(S_m(\tau_k),I_m(\tau_k))$. Hence the jump reproduces exactly the
prescribed biological impulse, without contamination from the continuous vector field.

\medskip

The hybrid system \eqref{eq:enhanced} is a system of Stieltjes differential equations with
\emph{two derivators} ---the identity $g(t)=t$ for the human compartments and $g_m$ for the mosquito
compartments. Its well-posedness follows from the theory of systems of Stieltjes differential
equations with several derivators \cite{LopezPousoMarquez2019}, and the integral formulation and the
decomposition into continuous and impulsive parts carry over to each derivator separately.

More explicitly, a solution of \eqref{eq:enhanced} is understood in integral form as
\[
S(t)=S(0)+\int_0^t -\beta_{hm}B S(s)\frac{I_m(s)}{N}\,ds,
\]
\[
I(t)=I(0)+\int_0^t \left(\beta_{hm}B S(s)\frac{I_m(s)}{N}-\gamma I(s)\right)ds,
\qquad
R(t)=R(0)+\int_0^t \gamma I(s)\,ds,
\]
while
\[
S_m(t)=S_m(0)+\int_{[0,t)}\mathcal F^{S_m}(s,I(s),S_m(s),I_m(s))\,d\mu_{g_m}(s),
\]
\[
I_m(t)=I_m(0)+\int_{[0,t)}\mathcal F^{I_m}(s,I(s),S_m(s),I_m(s))\,d\mu_{g_m}(s).
\]
This makes explicit that the human equations are integrated with respect to chronological time, whereas the mosquito equations are integrated with respect to \(\mu_{g_m}\).

In the full Stieltjes digital twin two impulsive mechanisms are present. The update times \(t_s^j\) carry the data-assimilation corrections, which act on the full state through the correction maps \(\mathcal C_j\) introduced in Section~\ref{sec:firstmodel}. The times \(\tau_k\) carry the biological impulses, which act only on the mosquito components through the \emph{pulse maps}
\[
\mathcal P_k(\mathbf x)
=
\bigl(S,\ I,\ R,\ S_m+J_k^S(S_m,I_m),\ I_m+J_k^I(S_m,I_m)\bigr)^{\top}.
\]
The relative position of the two families of atoms has to be specified, and two conventions are natural.
\begin{enumerate}
\item[(I)] \emph{Interior atoms.} Every \(\tau_k\) lies strictly between two consecutive update times. The two atomic supports are then disjoint, and the pre-update state \(\mathbf x(t_s^j)\) already contains the effect of any pulse that occurred in \((t_s^{j-1},t_s^j)\).
\item[(II)] \emph{Coincident atoms.} Some \(\tau_k\) coincide with update times, \(\tau_k=t_s^j\). The biological impulse is applied first and the data correction second, so that
\[
\mathbf x(\tau_k^+)=\mathcal C_j\bigl(\mathbf x_k^{\flat}\bigr),
\qquad
\mathbf x_k^{\flat}:=\mathcal P_k\bigl(\mathbf x(\tau_k)\bigr).
\]
The intermediate state \(\mathbf x_k^{\flat}\) is not a value of the solution, which at the atom only takes the values \(\mathbf x(\tau_k)\) and \(\mathbf x(\tau_k^+)\); it is a Lipschitz function of the left value \(\mathbf x(\tau_k)\). It is the model state that is compared with the observation at \(\tau_k\) (the one-step forecast) and it is the argument of the correction map.
\end{enumerate}
Both conventions are covered by Theorem~\ref{thm:full-stieltjes-dt-wp}, and Proposition~\ref{prop:conventions} shows that (II) is the limit of (I) when the biological atom approaches the update time from the left. The monthly numerical experiment uses convention (II), because its atoms \(\tau_k\) are integer month indices.

\medskip

\begin{theorem}[Well-posedness of the full Stieltjes digital twin]
\label{thm:full-stieltjes-dt-wp}
Let \(T>0\). Let
\[
T_s=\{t_s^0,\dots,t_s^J\},
\qquad
0=t_s^0<t_s^1<\cdots<t_s^J\le T,
\]
be the finite set of data-assimilation times, let \(T_s^\ast:=\{t_s^j:\ 1\le j\le J\}\) be the set of update times, excluding the initial time, and let
\[
T_m=\{\tau_k:\ k\in\mathcal K_m\}\subset(0,T]
\]
be the finite set of pairwise distinct biological mosquito impulse times. No relation between \(T_s^\ast\) and \(T_m\) is assumed. Put \(\mathcal T:=T_s^\ast\cup T_m\) and, for \(a\in\mathcal T\), define the jump map
\[
\Gamma_a
=
\begin{cases}
\mathcal P_k, & a=\tau_k\notin T_s^\ast,\\[2pt]
\mathcal C_j, & a=t_s^j\notin T_m,\\[2pt]
\mathcal C_j\circ\mathcal P_k, & a=t_s^j=\tau_k.
\end{cases}
\]
Let \(\alpha_j>0\), \(j=1,\dots,J\), and \(\delta_k>0\), \(k\in\mathcal K_m\). Define the human derivator
\[
g_H(t)=t+\sum_{\substack{1\le j\le J\\ t_s^j<t}}\alpha_j,
\]
and the mosquito derivator
\[
g_M(t)=t+\sum_{\substack{1\le j\le J\\ t_s^j<t}}\alpha_j+\sum_{\substack{k\in\mathcal K_m\\ \tau_k<t}}\delta_k,
\]
and denote by \(w_a:=\Delta g_M(a)\) the weight of the atom \(a\in\mathcal T\) of \(g_M\); thus \(w_a\) equals \(\alpha_j\), \(\delta_k\) or \(\alpha_j+\delta_k\) in the three cases above. The continuous part of both derivators is chronological time, which corresponds to the choice \(\alpha=1\) in \eqref{eq:derivator}. Assume that the epidemiological parameters are bounded and piecewise constant on the partition induced by \(\mathcal T\), that \(N>0\) and \(K>0\), that each correction map \(\mathcal C_j:\mathbb R^5\to\mathbb R^5\) is globally Lipschitz with \(\mathcal C_j(\mathcal A_N)\subset\mathcal A_N\), where \(\mathcal A_N=\{\mathbf x\in\mathbb R^5:\ S+I+R=N\}\), and that
the initial condition satisfies
\[
\mathbf x_0=(S_0,I_0,R_0,S_{m,0},I_{m,0})\in\mathbb R_+^5,
\qquad
S_0+I_0+R_0=N.
\]
Consider the system with several derivators
\begin{equation}\label{eq:full-MDT-several-derivators}
\begin{cases}
S'_{g_H}(t)=\mathcal F^S(t,\mathbf x(t)),\\
I'_{g_H}(t)=\mathcal F^I(t,\mathbf x(t)),\\
R'_{g_H}(t)=\mathcal F^R(t,\mathbf x(t)),\\
(S_m)'_{g_M}(t)=\mathcal F^{S_m}(t,\mathbf x(t)),\\
(I_m)'_{g_M}(t)=\mathcal F^{I_m}(t,\mathbf x(t)),
\end{cases}
\end{equation}
with \(\mathbf x=(S,I,R,S_m,I_m)\) and \(\mathbf x(0)=\mathbf x_0\). For \(t\notin T_s^\ast\) the human fields are given by the continuous malaria dynamics,
\[
\mathcal F^S(t,\mathbf x)
=
-\beta_{hm}B S\frac{I_m}{N},
\qquad
\mathcal F^I(t,\mathbf x)
=
\beta_{hm}B S\frac{I_m}{N}-\gamma I,
\qquad
\mathcal F^R(t,\mathbf x)=\gamma I,
\]
and for \(t\notin\mathcal T\) the mosquito fields are, with \(N_m=S_m+I_m\),
\[
\mathcal F^{S_m}(t,\mathbf x)
=
b\left(1+s_A\sin\left(\frac{2\pi t}{T_{\mathrm{season}}}+\phi\right)\right)
\left(1-\frac{N_m}{K}\right)N_m
-\beta_{mh}B S_m\frac{I}{N}
-\mu S_m,
\]
\[
\mathcal F^{I_m}(t,\mathbf x)
=
\beta_{mh}B S_m\frac{I}{N}
-
\mu I_m.
\]
At the atoms the fields are defined through the jump maps. Writing
\[
\Gamma_a=\bigl(\Gamma_a^{S},\Gamma_a^{I},\Gamma_a^{R},\Gamma_a^{S_m},\Gamma_a^{I_m}\bigr),
\]
we set
\[
\mathcal F^{X}(t_s^j,\mathbf x)=\frac{\Gamma_{t_s^j}^{X}(\mathbf x)-X}{\alpha_j},
\qquad X\in\{S,I,R\},\quad j=1,\dots,J,
\]
and
\[
\mathcal F^{X}(a,\mathbf x)=\frac{\Gamma_{a}^{X}(\mathbf x)-X}{w_a},
\qquad X\in\{S_m,I_m\},\quad a\in\mathcal T.
\]
In particular, at a biological atom \(\tau_k\notin T_s^\ast\) the mosquito fields are \(J_k^S(S_m,I_m)/\delta_k\) and \(J_k^I(S_m,I_m)/\delta_k\), as in \eqref{eq:mosquito-field-S}--\eqref{eq:mosquito-field-I}.
The solution is sought in the product space
\[
(S,I,R)\in \mathcal{AC}_{g_H}([0,T])^3,
\qquad
(S_m,I_m)\in \mathcal{AC}_{g_M}([0,T])^2,
\]
which is the natural class for a system with componentwise derivators.
Then there exists \(\tau\in(0,T]\) such that the full Stieltjes digital twin has
a unique solution on \([0,\tau]\) in the sense of systems of Stieltjes
differential equations with several derivators. This solution remains in \(\mathcal A_N\) and satisfies
\[
\mathbf x(a^+)=\Gamma_a\bigl(\mathbf x(a)\bigr)
\qquad\text{for every } a\in\mathcal T\cap[0,\tau).
\]
Moreover, if the solution remains bounded, then it can be continued
up to the final time \(T\).
\end{theorem}

\begin{proof}
We verify the hypotheses of the local existence and uniqueness theorem for
systems of Stieltjes differential equations with several derivators
\cite{LopezPousoMarquez2019}.

First, \(g_H\) and \(g_M\) are derivators. Indeed, both are nondecreasing and
left-continuous, and their associated Lebesgue--Stieltjes measures are
\[
\mu_{g_H}
=
\mu_{Id}
+
\sum_{j=1}^J \alpha_j\delta_{t_s^j},
\qquad
\mu_{g_M}
=
\mu_{Id}
+
\sum_{a\in\mathcal T}w_a\,\delta_{a}.
\]
Since \(\mathcal T\) is finite, both measures are finite on \([0,T]\). The initial time \(t_s^0=0\) is not an atom of either derivator.

Second, for every fixed state \(\mathbf x\), each component field is Borel measurable in \(t\), hence measurable
with respect to the corresponding Stieltjes measure: outside the finite atomic
sets the fields are piecewise continuous in time, because the parameters are
piecewise constant and the seasonal factor is continuous.

Third, the fields are Lipschitz with respect to the state variable on every ball
\(B_M=\{\mathbf x\in\mathbb R^5:\ \|\mathbf x\|\le M\}\). Indeed, since \(N\) is a fixed positive constant (Assumption~\ref{ass:model}(A2)), the continuous parts are polynomial in \(\mathbf x\) with bounded coefficients, and are therefore Lipschitz on \(B_M\), uniformly in \(t\). At the atoms, the pulse maps \(\mathcal P_k\) are affine and the correction maps \(\mathcal C_j\) are globally Lipschitz by hypothesis, so every jump map \(\Gamma_a\), including the compositions \(\mathcal C_j\circ\mathcal P_k\), is globally Lipschitz. Since the number of atoms is finite, the
corresponding Lipschitz constants are bounded, hence belong to \(L^1_{g_H}\) for the human components and to
\(L^1_{g_M}\) for the mosquito components. At a biological atom \(\tau_k\notin T_s^\ast\), only the mosquito equations have an atomic contribution; the human equations retain their chronological-time dynamics because \(\tau_k\) is not an atom of \(g_H\).

Fourth, the fields are bounded on \([0,T]\times B_M\): on the continuous part this follows
from the boundedness of the parameters, and at the finitely many atoms from the continuity of the jump maps and the positivity of the weights \(\alpha_j\) and \(w_a\). The measures being finite, these bounds belong to \(L^1_{g_H}\) and \(L^1_{g_M}\), respectively.

Therefore all hypotheses of the existence and uniqueness theorem for systems
with several derivators are satisfied. Hence there exists \(\tau\in(0,T]\) and
a unique solution on \([0,\tau]\).

The jump identity follows from the integral formulation: at an atom \(a\in\mathcal T\), each component changes by the value of its field multiplied by the weight of \(a\) in its own derivator. For the mosquito components this gives \(\Gamma_a^{X}(\mathbf x(a))-X(a)\). For the human components it gives the same expression if \(a\in T_s^\ast\), and zero if \(a=\tau_k\notin T_s^\ast\), which is again consistent because \(\mathcal P_k\) does not modify the human components. The solution stays in \(\mathcal A_N\) because the three human fields add up to zero away from the atoms, \(\mathcal C_j(\mathcal A_N)\subset\mathcal A_N\), and \(\mathcal P_k\) leaves \(S\), \(I\), \(R\) unchanged.

Finally, if the solution remains bounded, say in \(B_M\),
the same estimates remain valid with constants depending only on \(M\), and the local existence argument can be
restarted. Since \(T<\infty\) and the derivators have only finitely many atoms,
the solution can be continued step by step up to \(T\).
\end{proof}

\begin{remark}
Under convention (I) every jump map is either a pulse map \(\mathcal P_k\) or a correction map \(\mathcal C_j\). Under convention (II), at a coincident atom the jump map is \(\mathcal C_j\circ\mathcal P_k\), and the one-step forecast recorded by the implementation at that date is the intermediate state \(\mathbf x_k^{\flat}=\mathcal P_k(\mathbf x(\tau_k))\). The weights \(\alpha_j\) and \(\delta_k\) only normalise the atomic fields and do not affect the solution.
\end{remark}

\begin{proposition}[Consistency of the two conventions]
\label{prop:conventions}
Let the hypotheses of Theorem~\ref{thm:full-stieltjes-dt-wp} hold. Fix \(j\ge1\), write \(t_-=t_s^{j-1}\) and \(t_+=t_s^{j}\), and assume that \((t_-,t_+)\) contains no biological atom. Let \(\mathcal P\) be a pulse map with Lipschitz constant \(L_{\mathcal P}\), and let \(\mathbf y\) be the solution of the continuous system on \((t_-,t_+]\) issued from a given post-update state at \(t_-\). For \(0<\varepsilon<t_+-t_-\), let \(\mathbf x_\varepsilon\) be the solution issued from the same state when a single biological atom with pulse map \(\mathcal P\) is placed at \(t_+-\varepsilon\). Assume that \(\mathbf y\) and all \(\mathbf x_\varepsilon\) take their values in a compact set \(\Omega\subset\mathbb R^5\), and let \(M_\Omega\) be a bound on \((t_-,t_+]\times\Omega\) for the norm of the continuous vector field. Then
\[
\bigl\|\mathbf x_\varepsilon(t_+)-\mathcal P\bigl(\mathbf y(t_+)\bigr)\bigr\|
\le (1+L_{\mathcal P})\,M_\Omega\,\varepsilon .
\]
Consequently, as \(\varepsilon\to0^+\), the pre-update state of convention (I) converges to the intermediate state \(\mathbf x^{\flat}=\mathcal P(\mathbf y(t_+))\) of convention (II), and, by continuity of the correction map, the corresponding post-update states converge as well.
\end{proposition}

\begin{proof}
By uniqueness, \(\mathbf x_\varepsilon=\mathbf y\) on \((t_-,t_+-\varepsilon]\), and \(\mathbf x_\varepsilon((t_+-\varepsilon)^+)=\mathcal P(\mathbf y(t_+-\varepsilon))\). On \((t_+-\varepsilon,t_+]\) both \(\mathbf x_\varepsilon\) and \(\mathbf y\) solve the continuous system, whose field is bounded by \(M_\Omega\) on \(\Omega\); hence
\[
\bigl\|\mathbf x_\varepsilon(t_+)-\mathcal P(\mathbf y(t_+-\varepsilon))\bigr\|\le M_\Omega\,\varepsilon,
\qquad
\bigl\|\mathbf y(t_+)-\mathbf y(t_+-\varepsilon)\bigr\|\le M_\Omega\,\varepsilon .
\]
The Lipschitz continuity of \(\mathcal P\) and the triangle inequality give the estimate.
\end{proof}

The compactness assumption of Proposition~\ref{prop:conventions} holds, with \(\Omega=\mathcal H_N^K\), under the hypotheses of either of the next two results.

\medskip

\begin{corollary}[Global biologically admissible full Stieltjes digital twin]
\label{cor:full-stieltjes-dt-global}
Let
\[
\mathcal H_N^K=
\left\{
\mathbf x\in\mathbb R_+^5:
S+I+R=N,
\quad S_m+I_m\le K
\right\}.
\]
In addition to the hypotheses of Theorem \ref{thm:full-stieltjes-dt-wp}, assume that \(\mathbf x_0\in\mathcal H_N^K\), that \(0\le s_A\le1\), and that all correction and pulse maps preserve the admissible region,
\[
\mathcal C_j(\mathcal H_N^K)\subset\mathcal H_N^K,\quad j=1,\dots,J,
\qquad
\mathcal P_k(\mathcal H_N^K)\subset\mathcal H_N^K,\quad k\in\mathcal K_m.
\]
For the pulse maps \eqref{eq:mosquito-impulses} this holds, for instance, if
\[
\rho_k^S\ge -1,
\qquad
\rho_k^I\ge -1,
\qquad
\sigma_k\ge0,
\]
and
\[
(1+\rho_k^S)S_m+\sigma_k+(1+\rho_k^I)I_m\le K
\quad\text{whenever}\quad S_m,I_m\ge0,\ \ S_m+I_m\le K.
\]
Then the local solution given by Theorem \ref{thm:full-stieltjes-dt-wp} extends to the whole interval \([0,T]\) and remains in \(\mathcal H_N^K\). Under convention (II) the intermediate states \(\mathbf x_k^{\flat}\) also belong to \(\mathcal H_N^K\).
\end{corollary}

\begin{proof}
Between two consecutive atoms of \(g_M\), the system reduces to the continuous malaria equations with seasonal mosquito recruitment. The human total satisfies
\[
\frac{d}{dt}(S+I+R)=0,
\]
and the mosquito total satisfies
\[
\frac{dN_m}{dt}
=
\left[
 b\left(1+s_A\sin\left(\frac{2\pi t}{T_{\mathrm{season}}}+\phi\right)\right)
 \left(1-\frac{N_m}{K}\right)-\mu
\right]N_m.
\]
Since \(0\le s_A\le1\), the recruitment factor is nonnegative. Hence \([0,K]\) is positively invariant for the continuous mosquito-total equation: the vector field vanishes at \(N_m=0\) and is equal to \(-\mu K<0\) at \(N_m=K\). The quasi-positivity argument of Lemma~\ref{lem:positivity} gives the nonnegativity of the individual compartments during the continuous evolution.

At every atom \(a\) the state is mapped to \(\Gamma_a(\mathbf x(a))\), where \(\Gamma_a\) is a pulse map, a correction map or the composition \(\mathcal C_j\circ\mathcal P_k\); in all three cases \(\Gamma_a(\mathcal H_N^K)\subset\mathcal H_N^K\) by hypothesis, and in the third case \(\mathbf x_k^{\flat}=\mathcal P_k(\mathbf x(a))\in\mathcal H_N^K\) as well. Therefore the solution remains in the compact set \(\mathcal H_N^K\) after every continuous segment and after every atom, and the continuation criterion in Theorem \ref{thm:full-stieltjes-dt-wp} gives existence on the whole interval \([0,T]\).
\end{proof}

The uniform condition of Corollary~\ref{cor:full-stieltjes-dt-global} on the pulse maps fails for a purely additive emergence pulse (\(\rho_k^S=\rho_k^I=0\), \(\sigma_k>0\)), because a state with \(S_m+I_m=K\) is mapped outside the carrying-capacity region. Additive pulses are nevertheless admissible when mosquito mortality dominates recruitment near the carrying capacity and the pulses are sufficiently separated from the preceding atom, as the next result shows.

\begin{proposition}[Additive emergence pulses]
\label{prop:additive}
In addition to the hypotheses of Theorem \ref{thm:full-stieltjes-dt-wp}, assume that \(\mathbf x_0\in\mathcal H_N^K\), that \(0\le s_A\le1\), that \(\mathcal C_j(\mathcal H_N^K)\subset\mathcal H_N^K\) for \(j=1,\dots,J\), and that the pulses are additive,
\[
\rho_k^S=\rho_k^I=0,
\qquad
0\le\sigma_k\le\bar\sigma<K,
\qquad k\in\mathcal K_m.
\]
Let \(\bar b\) be an upper bound for the recruitment coefficient \(b\,(1+s_A)\) and \(\underline\mu>0\) a lower bound for the mortality \(\mu\) over all parameter intervals, and suppose that
\[
\kappa:=\underline\mu-\frac{\bar b\,\bar\sigma}{K}>0 .
\]
Set
\[
t^\ast:=\frac1\kappa\,\ln\frac{K}{K-\bar\sigma},
\]
and assume that every biological atom satisfies \(\tau_k-a_k\ge t^\ast\), where \(a_k:=\max\bigl((\{0\}\cup\mathcal T)\cap[0,\tau_k)\bigr)\) is the last atom preceding \(\tau_k\) (or the initial time). Then the solution exists on the whole interval \([0,T]\), remains in \(\mathcal H_N^K\), and satisfies
\[
N_m(\tau_k)\le K-\bar\sigma,
\qquad\text{hence}\qquad
N_m(\tau_k)+\sigma_k\le K,
\qquad k\in\mathcal K_m .
\]
\end{proposition}

\begin{proof}
We argue by induction over the atoms. Let \((a,a']\) be an interval between two consecutive elements of \(\{0\}\cup\mathcal T\), and assume that the state immediately after \(a\) belongs to \(\mathcal H_N^K\). As in the proof of Corollary~\ref{cor:full-stieltjes-dt-global}, the continuous flow keeps the solution in \(\mathcal H_N^K\) on \((a,a']\); in particular \(0\le N_m\le K\). Whenever \(N_m(t)\ge K-\bar\sigma\) we have \(0\le 1-N_m(t)/K\le\bar\sigma/K\), and therefore
\[
N_m'(t)
\le
\left(\frac{\bar b\,\bar\sigma}{K}-\underline\mu\right)N_m(t)
=-\kappa\,N_m(t)<0 .
\]
Two consequences follow. First, the set \(\{N_m\le K-\bar\sigma\}\) is positively invariant for the continuous flow. Second, as long as \(N_m\ge K-\bar\sigma\), we have \(N_m(t)\le N_m(a^+)\,e^{-\kappa(t-a)}\le K e^{-\kappa(t-a)}\). Since \(K e^{-\kappa t^\ast}=K-\bar\sigma\), we conclude that \(N_m(t)\le K-\bar\sigma\) for every \(t\in(a,a']\) with \(t-a\ge t^\ast\).

If \(a'=\tau_k\) is a biological atom, then \(a=a_k\) and \(\tau_k-a_k\ge t^\ast\), so \(N_m(\tau_k)\le K-\bar\sigma\). The additive pulse then yields nonnegative mosquito compartments with total \(N_m(\tau_k)+\sigma_k\le K\), that is, \(\mathcal P_k(\mathbf x(\tau_k))\in\mathcal H_N^K\). If \(a'\) is also an update time, the correction map is applied to this state and returns an element of \(\mathcal H_N^K\) by hypothesis; the same holds if \(a'\) is an update time only. Thus the state immediately after \(a'\) belongs to \(\mathcal H_N^K\), which closes the induction. The solution therefore remains in the compact set \(\mathcal H_N^K\), and the continuation criterion of Theorem \ref{thm:full-stieltjes-dt-wp} gives existence on \([0,T]\).
\end{proof}

\medskip

In the numerical implementation, the seasonal forcing term is constrained to remain nonnegative. The continuous mosquito dynamics are advanced with respect to the continuous part of \(g_m\). At an atom \(\tau_k\), the continuous vector field is not multiplied by the jump \(\delta_k\); instead, the prescribed impulse is applied directly through
\[
S_m(\tau_k^+)-S_m(\tau_k)=J_k^S(S_m(\tau_k),I_m(\tau_k)),
\qquad
I_m(\tau_k^+)-I_m(\tau_k)=J_k^I(S_m(\tau_k),I_m(\tau_k)),
\]
which is equivalent to assigning the value \(J_k/\delta_k\) to the Stieltjes vector field at the atom. As a safeguard, the code raises an error if an impulse produces a negative compartment or leaves the region \(S_m+I_m\le K\).

\begin{table}[ht!]
\centering
\caption{Numerical settings used by the monthly Python experiment. Bounds apply only to rolling-window fitting.}
\label{tab:stieltjes-params}
\scriptsize
\begin{tabular}{|c|c|p{6.4cm}|}
\hline
\textbf{Quantity} & \textbf{Value or bounds} & \textbf{Role} \\
\hline
$N$ & $4{,}377{,}770$ & Fixed human population \\
$K$ & $8.0\times10^7$ & Mosquito carrying capacity \\
$B$ & $0.5\,\mathrm{month}^{-1}$ & Fixed nominal biting-rate scale; only $\beta_{hm}B$ and $\beta_{mh}B$ are identifiable \\
$\mu$ & $2.14\,\mathrm{month}^{-1}$ & Fixed mosquito mortality \\
$T_{\rm season}$ & $12$ months & Annual mosquito recruitment period \\
$s_A$ & $0.40$ & Seasonal amplitude; nonnegative forcing \\
$\phi$ & data-derived & Places the synthetic mosquito peak about one month before the case peak \\
$\tau_k$ & annual, data-derived month indices & Rainy-season onset atoms; they coincide with update times (convention (II)) \\
$\delta_k$ & $0.5$ & Derivator jump; the atomic field is $J_k/\delta_k$ \\
$\sigma$ & $3.0\times10^6$ mosquitoes & Additive susceptible-mosquito emergence pulse \\
$\rho^S,\rho^I$ & $0,0$ & No proportional impulse in the reported run \\
$q$ & $1$ (configurable) & Observation proxy $Y_k=qI(t_s^k)+\varepsilon_k$ \\
\hline
$\beta_{hm}$ & initial $0.02$; $[10^{-3},3]$ & Fitted effective mosquito-to-human transmission coefficient \\
$\beta_{mh}$ & initial $0.03$; $[10^{-4},1]$ & Fitted effective human-to-mosquito transmission coefficient \\
$\gamma$ & initial $0.25$; $[0.05,1]$ month$^{-1}$ & Fitted human recovery rate \\
$b$ & initial $4.0$; $[2.2,10]$ month$^{-1}$ & Fitted mosquito recruitment coefficient \\
\hline
\end{tabular}
\end{table}

Table~\ref{tab:parameters} gives only day-based literature ranges for context, whereas Table \ref{tab:stieltjes-params} lists the values used by the executed code, which is monthly throughout. The reported run uses purely additive pulses, so Corollary~\ref{cor:full-stieltjes-dt-global} does not apply to it, but Proposition~\ref{prop:additive} does. Indeed, with $\bar\sigma=3\times10^6$, $K=8\times10^7$, $\underline\mu=2.14$ and $\bar b=10\,(1+s_A)=14$ (the upper fitting bound of $b$ times the maximal seasonal factor) one obtains
\[
\kappa=2.14-\frac{14\times 3\times10^{6}}{8\times10^{7}}=1.615\ \mathrm{month}^{-1},
\qquad
t^\ast=\frac{1}{\kappa}\ln\frac{K}{K-\bar\sigma}\approx 0.024\ \text{months},
\]
whereas every pulse is applied one month after the preceding update. Moreover, the correction maps preserve $\mathcal H_N^K$, because the observed case counts are far below $N$ and the synthetic mosquito observations stay below $K$ (Figure~\ref{fig:mosquito-fit}). Hence the solution remains in $\mathcal H_N^K$ and the safeguard of the code cannot be triggered, in agreement with the executed run.

\subsection{Numerical Implementation}\label{sec:numerics}

The Stieltjes twin uses the same monthly rolling-calibration code as the classical twin. The only structural changes are the annual seasonal factor in mosquito recruitment and the biological impulses, applied according to convention (II). Both runs share the same human observations, synthetic mosquito observations, fitting windows, parameter bounds, optimizer, state corrections, ensemble size and evaluation months. This comparison holds the computational protocol fixed but does not isolate the effect of the Stieltjes representation, since the model with Stieltjes dynamics additionally includes seasonal recruitment and impulses.

\emph{Synthetic entomological recovery experiment.} Because no measured Kebbi
\emph{Anopheles} series is available, the target vector abundance is generated from
\[
N_m'=b_{\rm true}s(t)\left(1-\frac{N_m}{K}\right)N_m-\mu N_m,
\qquad
s(t)=\max\left\{0,1+s_A\sin\left(\frac{2\pi t}{12}+\phi\right)\right\},
\]
with $b_{\rm true}=4.5\,\mathrm{month}^{-1}$, annual additive pulses $N_m(\tau_k^+)=N_m(\tau_k)+3\times10^6$, and multiplicative log-normal observation noise with coefficient of variation $8\%$. A seasonal sporozoite fraction splits $N_m$ into $S_m$ and $I_m$. Since the infection-transfer terms cancel in the equation for $N_m=S_m+I_m$, this scalar generator is consistent with the total-mosquito balance of the five-compartment model.

Two features of this experiment determine how its results must be read. First, the seasonal phase $\phi$ and the pulse months $\tau_k$ of the generator are obtained from the monthly climatology of the complete case series (the mosquito peak is placed one month before the climatological case peak, and the pulses two months before the mosquito peak), and the same $\phi$ and $\tau_k$ are supplied to the Stieltjes twin. The experiment therefore tests the recovery of the recruitment amplitude, of the transmission coefficients and of the pulse magnitude when the seasonal phase and the pulse schedule are known; it does not test the identification of the schedule. Second, the generator depends on the case series only through $\phi$ and $\tau_k$, and these coincide with the values obtained from the first twelve months (January--December 2015) alone. A re-run in which the generator is frozen on that first year reproduces every synthetic observation, fitted trajectory and forecast to rounding error. Hence, for the evaluation months from January 2016 onwards, the synthetic target and the forecasts built on it are non-anticipative; restricting Table~\ref{tab:metrics} to these months changes every positive $R^2$ value by at most $0.04$, lowers further the already negative values of the mechanistic-only human forecasts, and preserves every ordering reported below. The mosquito ``one-step forecast'' in Table~\ref{tab:metrics} nevertheless remains a one-step recovery with known phase and schedule.

\subsection{Results and Discussion}\label{sec:results}

Table~\ref{tab:metrics} collects the performance of both twins on the common post-baseline evaluation period.

\begin{figure}[ht!]
\centering
\includegraphics[width=\textwidth]{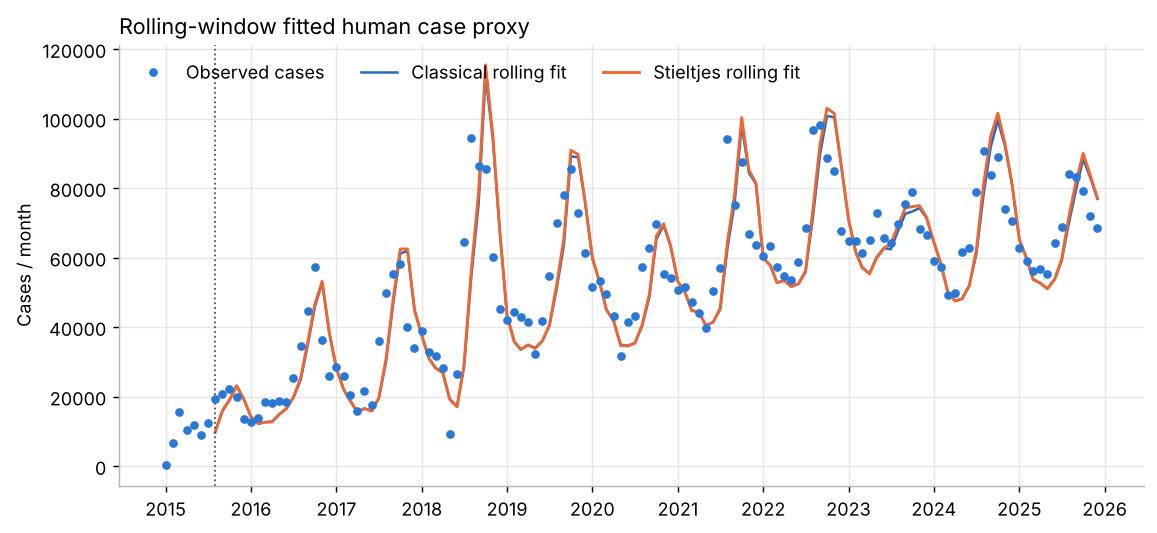}
\caption{Rolling-window fitted human-case proxy for the classical and Stieltjes digital twins. Blue points are the 132 observed monthly confirmed-case counts. The dotted vertical line marks the end of the seven-month baseline. The two fitted trajectories are nearly indistinguishable in aggregate performance.}
\label{fig:results-classical}
\end{figure}

\begin{figure}[ht!]
\centering
\includegraphics[width=\textwidth]{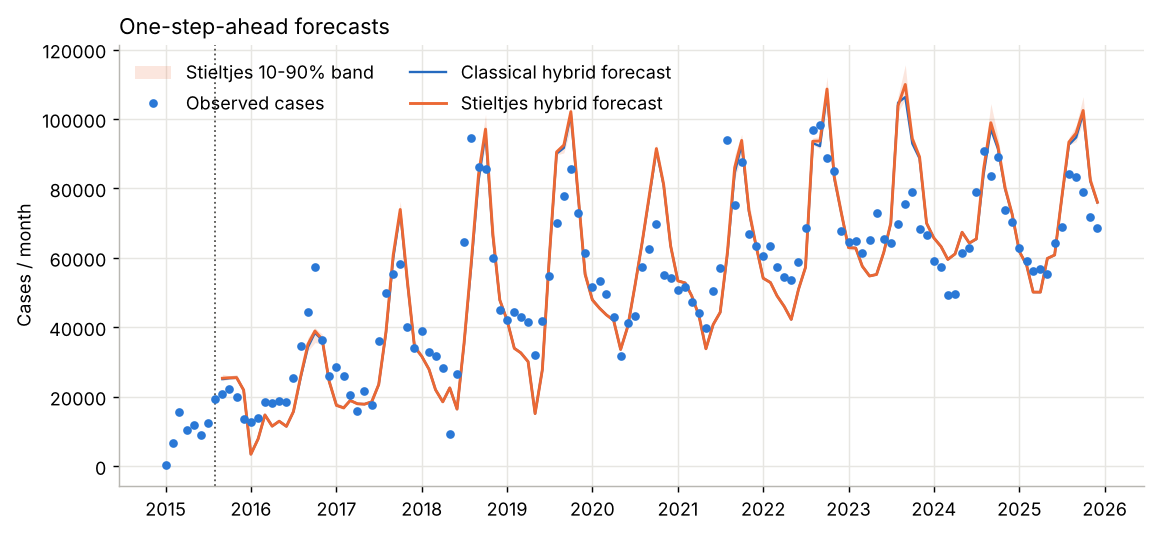}
\caption{One-month-ahead human forecasts. The hybrid forecast combines the mechanistic ensemble with a past-only lag-12 seasonal-naive term; the mechanistic component is corrected to the synthetic mosquito series described in Section~\ref{sec:numerics}. The shaded region is the Stieltjes $10$--$90\%$ parameter-perturbation envelope, not a calibrated predictive interval.}
\label{fig:forecast}
\end{figure}

\emph{Human proxy.} Figures \ref{fig:results-classical} and~\ref{fig:forecast} compare the two formulations. Their performances are very similar: the rolling-window $R^2$ values are $0.7095$ (classical) and $0.7131$ (Stieltjes), whereas the hybrid forecast gives $0.7218$ and $0.7136$, respectively, so the classical hybrid forecast is marginally better in this run. The mechanistic-only human forecast has negative $R^2$ in both formulations. Table~\ref{tab:metrics} also reports two reference predictors evaluated on the same months: the seasonal-naive forecast (the value observed twelve months earlier), with $R^2=0.6226$, and last-month persistence (the most recent observation), with $R^2=0.7616$. Both hybrid forecasts improve on the seasonal-naive predictor, but neither improves on persistence. The improvement over the seasonal-naive predictor is explained by the fact that the mechanistic forecast is issued from the state corrected to the current observation and therefore carries the current level: its correlation with the most recent observation is $0.84$ (classical) and $0.83$ (Stieltjes). The human forecast skill thus reflects the information contained in the last observations rather than the transmission mechanism. These results do not establish a human-case advantage of the Stieltjes model, and no such advantage is claimed.

\begin{table}[ht!]
\centering
\caption{Performance on the common post-baseline evaluation period. Human errors are in cases/month and vector errors in mosquitoes. ``Rolling fit'' is a window-fitted trajectory and is not an independent forecast. The mosquito ``one-step forecast'' is a one-step recovery of a synthetic series whose seasonal phase and pulse schedule are known to the model (Section~\ref{sec:numerics}). The reference rows are the observation twelve months earlier (the most recent observation when no twelve-month lag is available) and last-month persistence, evaluated on the same months.}
\label{tab:metrics}
\scriptsize
\begin{tabular}{|l|l|l|r|r|r|}
\hline
\textbf{Model} & \textbf{Target} & \textbf{Mode} & \textbf{RMSE} & \textbf{MAE} & \textbf{$R^2$} \\
\hline
Classical & Human & rolling fit & $1.1606\times10^4$ & $8.6472\times10^3$ & $0.7095$ \\
Classical & Human & mechanistic forecast & $3.9615\times10^4$ & $2.5147\times10^4$ & $-2.3846$ \\
Classical & Human & hybrid forecast & $1.1359\times10^4$ & $8.7845\times10^3$ & $0.7218$ \\
\hline
Stieltjes & Human & rolling fit & $1.1533\times10^4$ & $8.5703\times10^3$ & $0.7131$ \\
Stieltjes & Human & mechanistic forecast & $4.1733\times10^4$ & $2.6288\times10^4$ & $-2.7562$ \\
Stieltjes & Human & hybrid forecast & $1.1524\times10^4$ & $8.8965\times10^3$ & $0.7136$ \\
\hline
Reference & Human & seasonal naive (lag 12) & $1.3229\times10^4$ & $1.0434\times10^4$ & $0.6226$ \\
Reference & Human & persistence (lag 1) & $1.0515\times10^4$ & $7.6614\times10^3$ & $0.7616$ \\
\hline
Classical & Mosquito total & rolling fit & $9.7285\times10^6$ & $8.0888\times10^6$ & $0.3454$ \\
Classical & Mosquito total & one-step forecast & $1.2742\times10^7$ & $1.0772\times10^7$ & $-0.1229$ \\
Stieltjes & Mosquito total & rolling fit & $4.0234\times10^6$ & $3.1678\times10^6$ & $0.8880$ \\
Stieltjes & Mosquito total & one-step forecast & $4.7248\times10^6$ & $3.7188\times10^6$ & $0.8456$ \\
\hline
Reference & Mosquito total & seasonal naive (lag 12) & $5.5770\times10^6$ & $4.0478\times10^6$ & $0.7849$ \\
Reference & Mosquito total & persistence (lag 1) & $8.0151\times10^6$ & $6.2611\times10^6$ & $0.5557$ \\
\hline
\end{tabular}
\end{table}

\begin{figure}[ht!]
\centering
\includegraphics[width=\textwidth]{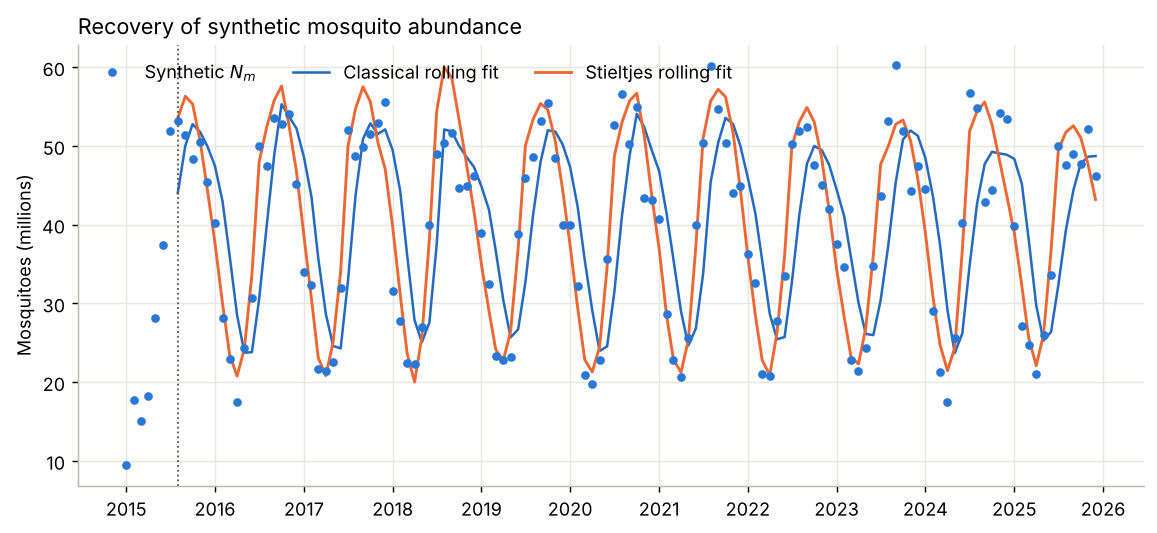}
\caption{Recovery of the synthetic mosquito abundance by the classical and Stieltjes rolling fits. The Stieltjes formulation reproduces the annual pulse structure used by the synthetic generator.}
\label{fig:mosquito-fit}
\end{figure}

\emph{Synthetic vector compartment.} The central comparative result occurs in the synthetic vector compartment. The classical model obtains $R^2=0.345$ for the rolling fit and $R^2=-0.123$ for the one-step recovery. The Stieltjes model, which contains the seasonal recruitment and the pulse mechanism used by the generator, obtains $R^2=0.888$ and $R^2=0.846$, respectively; for comparison, the seasonal-naive predictor attains $R^2=0.785$ and last-month persistence $R^2=0.556$ on the synthetic series. Figure~\ref{fig:mosquito-fit} shows the corresponding trajectories. Two qualifications are essential. First, the classical twin of Section~\ref{sec:firstmodel} contains neither seasonal recruitment nor impulses, so the comparison quantifies the value of including these two mechanisms, and it does not apportion the gain between them; a seasonal-only ablation is required for that purpose. Second, by the equivalence between impulsive systems and Stieltjes equations established in Section~\ref{sec:firstmodel}, a classical model endowed with the same seasonal recruitment and the same impulses, integrated as an impulsive system, produces exactly the same trajectories as the Stieltjes twin; the advantage reported here is therefore a property of the modelled mechanisms, which the Stieltjes derivator represents exactly and unifies with the assimilation atoms, and not a numerical advantage of the derivative calculus itself. The experiment is a successful structural-recovery test, with known phase and schedule, for a data-generating mechanism that contains pulses; it is not an external validation against real entomological measurements, and it is not evidence that the mosquito abundance in Kebbi State is discontinuous.

\begin{figure}[ht!]
\centering
\includegraphics[width=\textwidth]{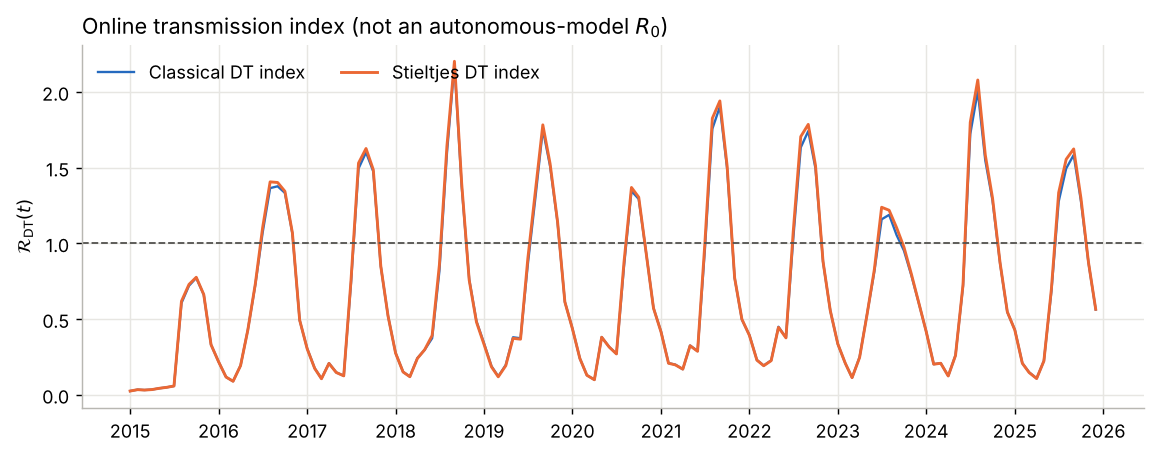}
\caption{Online digital-twin transmission indicator $\mathcal R_{\mathrm{DT}}(t)$ for both models. This current-state indicator is distinct from the autonomous basic reproduction number $\mathcal R_0$ in Proposition~\ref{prop:R0}.}
\label{fig:R0}
\end{figure}

Figure~\ref{fig:R0} displays the online indicator $\mathcal R_{\mathrm{DT}}$ defined in \eqref{eq:RDT}. The two runs give almost identical curves, because both are corrected to the same mosquito observations; the indicator ranges approximately from $0.025$ to $2.20$ and crosses the reference level one seasonally. Because current corrected mosquito abundance replaces the disease-free equilibrium abundance, the line $\mathcal R_{\mathrm{DT}}=1$ is used as a descriptive transmission-potential reference. It is not asserted to be a rigorous stability threshold for the nonautonomous assimilated system.

\begin{remark}
The Python implementation aggregates mosquito life stages into $S_m$ and $I_m$. The theoretical several-derivator formulation remains valid for this aggregation, but a stage-structured entomological validation would require additional measured states and parameters.
\end{remark}

\begin{remark}[Evaluation protocol]
The rolling fit is an in-sample (assimilation) measure: every observation in the window is used for calibration and state correction. The one-step forecasts are out of sample by one month, since the observation at $t_s^{k+1}$ is not used when the forecast is issued at $t_s^k$. No longer horizon is assessed and no hold-out period is reserved, because every observation is eventually assimilated. A multi-horizon evaluation on a withheld period is the natural next step and is left for future work.
\end{remark}

\section{Conclusions and Future Work}\label{sec:conclusions}

We have formulated classical and Stieltjes digital-twin versions of a mosquito--human malaria model and compared them with a common Python implementation. Three objects have been distinguished. The first is the autonomous invasion threshold $\mathcal R_0$ of the classical model at its fully susceptible disease-free equilibrium (Proposition~\ref{prop:R0}), together with its dependence on the chosen disease-free equilibrium through the factor $\sqrt{S^\ast/N}$. The second is the online indicator $\mathcal R_{\mathrm{DT}}(t)$, which evaluates the same expression at the current corrected mosquito abundance; it is a descriptive measure of transmission potential and not a stability threshold. The third is the repeatedly corrected Main Digital Twin, for which we have proved local well-posedness on the affine state space of constant human population and global well-posedness under positivity and invariant-region hypotheses, both when biological impulses fall between data updates and when they coincide with them.

\medskip
\textbf{Scope and limitations.} Neither $\mathcal R_0$ nor $\mathcal R_{\mathrm{DT}}(t)$ is a per-cycle growth multiplier of the impulsive system; such a quantity would require a periodic next-generation analysis in the Stieltjes setting, which is not undertaken here. Likewise, only the linear stability of the infected subsystem is asserted, and the global nonlinear dynamics of the model are not analysed. These are natural directions for future work.

The numerical study gives two different levels of evidence. For human data, confirmed monthly case counts are used only through the provisional proxy $Y_k=qI(t_s^k)$ with $q=1$. The classical and Stieltjes human results are essentially comparable, and the mechanistic-only forecasts are not competitive. The hybrid forecasts improve on a seasonal-naive predictor but not on last-month persistence, so the human forecast skill reflects the information in the most recent observations rather than the transmission mechanism. For mosquitoes, the observations are synthetic and contain seasonal recruitment and
annual emergence pulses. In this controlled recovery setting, with the phase and the schedule known to the model, the twin endowed with these two mechanisms clearly outperforms the unforced classical model. This supports the correctness of the implementation and the usefulness of representing seasonal and impulsive recruitment when such events are present, but it does not apportion the gain between the two mechanisms, it attributes no numerical advantage to the Stieltjes calculus itself, and it does not validate the mechanism for Kebbi State.

\begin{figure}[ht!]
\centering
\includegraphics[width=\textwidth]{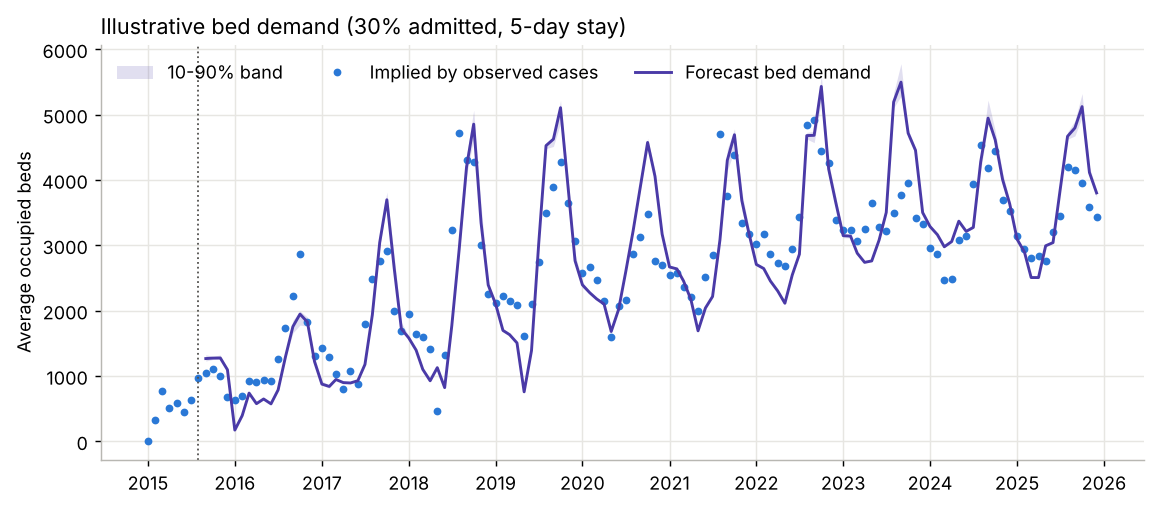}
\caption{Illustrative average bed demand obtained from the Stieltjes human forecast under the hypothetical scenario of $30\%$ admitted and five-day mean stay. The shaded region propagates the $10$--$90\%$ parameter-perturbation envelope; it is illustrative, not a calibrated predictive interval.}
\label{fig:hospital}
\end{figure}

The current forecast horizon is one month. Statements about predicting peaks several weeks or months in advance would require a separate multi-horizon, strictly out-of-sample study. Similarly, hospital demand is not observed in the dataset, and the case series consists of uncomplicated malaria, which is normally managed on an outpatient basis. Figure~\ref{fig:hospital} therefore only illustrates how a case forecast can be converted into a bed-occupancy curve, under a purely hypothetical scenario in which $30\%$ of monthly reported cases are admitted with a five-day average stay; this gives a forecast peak of approximately $5{,}500$ average occupied beds, with the upper ensemble curve below approximately $5{,}800$. These figures are not a capacity estimate for Kebbi State: a usable conversion requires observed admission fractions, lengths of stay, reporting coverage and delay distributions.

The next steps are to introduce an explicit incidence/reporting observation model, so that the mechanistic forecast can add information beyond persistence; to carry out a seasonal-only ablation that apportions the vector gain between seasonal recruitment and emergence pulses; to reserve a genuine hold-out period and evaluate several forecast horizons; to obtain entomological, climate and intervention data; to estimate, rather than prescribe, the biological impulse schedule; and to replace the illustrative hospital conversion by an observation-based resource model. A spatial or age-structured extension should be considered only after these identifiability and validation questions have been resolved.

\medskip
\noindent\textbf{Data availability.}
The monthly confirmed-malaria counts for Kebbi State used in this study (January 2015--December 2025) are available from the corresponding author upon reasonable request, subject to the permission of the data provider. The period January 2015--May 2024 corresponds to the case series analysed in \cite{Bakare2026}. The synthetic entomological series is generated by the procedure described in Section~\ref{sec:numerics}.

\medskip
\noindent\textbf{Code availability.}
The Python source code and the verification scripts are available from the corresponding author upon reasonable request.

\end{document}